\documentstyle[11pt]{article}

\newtheorem{theorem}{Theorem}[subsection]
\newtheorem{lemma}[theorem]{Lemma}
\newtheorem{corollary}[theorem]{Corollary}
\newtheorem{proposition}[theorem]{Proposition}

\newtheorem{rem}[theorem]{Remark}
        \newenvironment{remark}{\begin{rem}\rm}{\end{rem}}
\newtheorem{define}[theorem]{Definition}
        \newenvironment{definition}{\begin{define}\rm}{\end{define}}
\newtheorem{ex}[theorem]{Example}%
        \newenvironment{example}{\begin{ex}\rm}{\end{ex}}

\newcommand{\1}{_{(1)}}
\newcommand{\2}{_{(2)}}

\font\twlmsbm=msbm10 scaled \magstep1
\font\egtmsbm=msbm8
\font\sixmsbm=msbm6
\newfam\msbmfam

\textfont\msbmfam=\twlmsbm
\scriptfont\msbmfam=\egtmsbm
\scriptscriptfont\msbmfam=\sixmsbm

\newcommand{\calI}{{\cal I}}
\newcommand{\calS}{{\cal S}}
\newcommand{\calB}{{\cal B}}
\newcommand{\calH}{{\cal H}}
\newcommand{\Mdi}{M_{d_i}(\Bbb{C})}
\newcommand{\eps}{\varepsilon}
\newcommand{\Tr}{\mbox{Tr}\,}
\newcommand{\id}{\mbox{id}{}_M}

\newcommand{\rank}{\mbox{rank}}

\def\underset#1#2{\mathop{#2}\limits_{#1}}
\def\text#1{\hbox{\rm #1}}

\font\twlmsbm=msbm10 scaled \magstep1
\font\egtmsbm=msbm8
\font\sixmsbm=msbm6
\newfam\msbmfam

\textfont\msbmfam=\twlmsbm
\scriptfont\msbmfam=\egtmsbm
\scriptscriptfont\msbmfam=\sixmsbm

\def\Bbb#1{{\fam\msbmfam\relax#1}}

\title{\bf ALGEBRAIC VERSIONS OF A FINITE-DIMENSIONAL QUANTUM GROUPOID}

\author{Dmitri Nikshych\thanks{UCLA, Department of Mathematics, 
405 Hilgard Avenue, Los Angeles, CA
90095-1555;E-mail: nikshych@math.ucla.edu} \and Leonid Vainerman\thanks
{URA CNRS (Case 191), Universit\'e Pierre et Marie Curie, 4, place Jussieu,
F-75252, Paris Cedex 05 France; E-mail: vainerma@math.jussieu.fr}}

\date{October 28, 1998}

\begin{document}

\maketitle
\vskip 0.5cm
\begin{abstract}
We establish the equivalence of three approaches to the theory of finite
dimensional quantum groupoids. These are the generalized Kac algebras of
T.~Yamanouchi, the weak Kac algebras, i.e., the weak $C^*$-Hopf algebras
introduced by G. B\"ohm--F.~Nill--K.~Szlach\'anyi which have an involutive
antipode, and the Kac bimodules. The latter are an algebraic version of the Hopf
bimodules of J.-M.~Vallin.  We also study the structure and construct examples
of finite dimensional  quantum groupoids.
\end{abstract}

\bigskip\bigskip
\centerline{\bf INTRODUCTION}
\bigskip
 
{\bf 0.1.} Our starting point is the following theorem, which was
conjectured by A.~Ocneanu (see \cite{O}, the appendix in \cite{ES} and 
remark 4.7a in \cite{GHJ}) and was proved by
W.~Szymanski \cite{Szym}, R. Longo \cite{L}, and M.-C. David \cite{D}:
 
If $M_0\subset M_1\subset M_2\subset M_3$ is a {\it Jones's tower} of
type II${}_1$ von Neumann factors with finite index \cite{GHJ}, which is
{\it irreducible} (i.e., $M'_0
\cap M_1=C$) and {\it of depth 2} (i.e., $M'_0\cap M_3$ is a factor),
then $M'_1
\cap M_3$ has a natural finite-dimensional (f-d) Kac algebra \cite{KP},
\cite{ES} structure which acts {\it outerly} on $M_2$ in such a way
that the resulting crossed product is isomorphic to $M_3$. This gives an
intrinsic characterization of inclusions of the mentioned type.

A similar result for infinite index inclusions was conjectured in
\cite{HO} and proved for arbitrary type factors in \cite{EN}, \cite{E}.

It is natural to extend the above theorem to include
{\it reducible} (i.e., $\Bbb{C}\subset M'_0\cap M_1$) inclusions of 
depth 2. It is easy to produce an example of such an inclusion by
taking a crossed product $M_0\subset M_0\times_\alpha G$ with $G$  a 
finite group and  $\alpha$ an action which is not outer \cite{EN},
\cite{EVal}, \cite{NSzW}. In this more general context, $M'_1\cap M_3$ 
is no longer a Kac algebra. Nevertheless, it has an algebraic
structure which can be used in much the same manner. 
\medskip

{\bf 0.2.}  A f-d {\it generalized (gen.) Kac algebra} \cite{Y} is 
a $C^*$-bialgebra with a (generally) non-unital coproduct, an involutive 
antipode and a {\it Haar trace}. Exactly as for Kac algebras  \cite{KP}, 
\cite{ES}, the corresponding duality theory is based on the usage of the 
{\it fundamental operator} $W$. $W^*$ satisfies the {\it pentagonal
relation} \cite{BS}, but it is a partial isometry rather than
necessarily being a unitary. Despite this fact,
one can still construct a gen.\ Kac algebra from a given f-d
multiplicative partial isometry \cite{Val3}. 
Since a gen.\ Kac algebra is commutative iff it is
isomorphic to the algebra of complex-valued functions on a finite
groupoid (see \cite{Y}, Theorem 7.32), one can consider two above
notions as ``measurable versions'' of a f-d {\it quantum groupoid} --
``measurable'', because their definitions contain explicitly (a gen. 
Kac algebra) or implicitly (a multiplicative partial isometry) a Haar 
measure.

The corresponding infinite-dimensional structures are {\it Hopf
bimodules} \cite{Val1} and {\it pseudo-multiplicative unitaries}
\cite{Val2} introduced using a fibered product of von Neumann algebras
\cite{S}. In \cite{EVal}, a structure theory was considered
for von Neumann algebraic inclusions of depth 2 and any index,
equipped with an operator-valued weight verifying certain regularity conditions. 
On the other hand, a f-d {\it weak $C^*$-Hopf algebra} (\cite{BNSz},
\cite{BSz}) generalizing a Kac algebra and formulated in purely
algebraic terms was used in \cite{NSzW} in order to show
that any inclusion given by a crossed product with such an algebra is of
depth 2. These authors employed a  $C^*$-bialgebra with an antipode and 
a counit instead of a Haar measure. It turns out that any depth 2
finite index inclusion of II${}_1$ factors can be characterized this way - see
\cite{NV}, where we used in an essential way the results of the present paper
and those of \cite{N} on actions of weak Kac algebras (i.e., weak
$C^*$-Hopf algebras with an involutive antipode) on operator algebras.

We note that other formulations of the notion of quantum groupoids
have been proposed in \cite{Lu}, \cite{Mal}, and \cite{V2}.

In what follows all $C^*$-algebras and linear spaces over $\Bbb{C}$ are
f-d.
\medskip  

{\bf 0.3.} In this paper we establish the equivalence of three versions
of a f-d quantum groupoid: the gen.\ Kac algebra, the weak Kac algebra, 
and the Kac bimodule, an ``algebraization'' of the Hopf bimodule. 
We also study the structure and
construct examples of f-d quantum groupoids (some other examples are
given in \cite{EVal}, \cite{BSz}, \cite{BNSz}).

The paper is organized as follows. In Section 1 (Preliminaries) we discuss,
following \cite{S}, \cite{EVal}, \cite{Val3}, a relative tensor product of
Hilbert modules and a fibered product of f-d $C^*$-algebras. Next,
important notions of {\it Cartan subalgebras} (this term is borrowed
from \cite{R}, \cite{Val3}) and {\it counital maps} in $C^*$-bialgebras
are introduced. In their terms we define two auxiliary structures --
counital $C^*$-bialgebras and Hopf bimodules, study their elementary
properties, and prove their equivalence.

In Section 2 we first discuss weak Kac algebras (see also \cite{BNSz},
\cite{BSz} for much more complicated theory of general weak $C^*$-Hopf
algebras). The main results are existence and uniqueness of a Haar projection,
normalized Haar trace and a Haar conditional expectation on any weak Kac
algebra. This material can be read independently, one can look at Preliminaries
if necessary.
Then the equivalence of the categories of gen.\ and weak Kac algebras is
established (see also a remark in \cite{Nill} on this subject). The equivalent
language of Kac bimodules (i.e., Hopf bimodules \cite{EVal}, \cite{Val1}, 
\cite{Val3} equipped with a counit) gives a possibility to extend this theory 
to infinite-dimensional case.

It is natural to consider a f-d quantum groupoid as non-trivial if it is
non-commutative and non-cocommutative (i.e., does not come from groupoids)
and is not a usual Kac algebra. Such examples are presented in Section 3.
In particular, a finite group acting on a weak Kac algebra gives, via
crossed product construction, a new weak Kac algebra. In such a way, we
construct a series of concrete non-trivial quantum groupoids of
dimensions $n^3\ (n\geq 2)$.

Then we show that for any f-d $C^*$-algebra $A$ there exists a unique weak
Kac algebra structure (non-trivial when $A$ is non-abelian) on the full
matrix algebra $M_n(\Bbb{C})$ (where $n=\dim A$) with Cartan subalgebras
isomorphic to $A$. This gives a classification of weak Kac algebras simple
as $C^*$-algebras.

{\bf Acknowledgements.} 
The authors would like to thank E.~Effros for stimulating discussions
and E.~Vaysleb for his useful comments.
The second author is deeply grateful to M.~Enock and J.-M.~Vallin for many
important discussions and also to Max-Planck-Institut f\"{u}r Mathematik in Bonn
for excellent conditions for his work.


\begin{section} {Preliminaries}

\begin{subsection} {Relative tensor product of Hilbert modules}
\label{relative tensor}

{\bf 1.1.1.} Let $\theta$ be a faithful trace on a $C^*$-algebra $N$,
$L_2(N)$ a Hilbert space given by the GNS-construction for $(N,\theta)$,
$\lambda$ (resp., $\rho$) a faithful unital $*$-representation
(resp., anti-$*$-representation) of $N$ on $L_2(N)$ by operators of
left (resp., right) multiplication. Clearly, the commutant
$\lambda(N)'=\rho(N)$.

Let $\alpha$ be a faithful unital $*$-representation of $N$ in some
other Hilbert space $H$. Then $H$ can be viewed as a left $N$-module
$_\alpha H:n\cdot_\alpha\zeta=\alpha(n)\zeta\ (\forall n\in N,\zeta\in
H)$. An operator $R(\zeta):L_2(N)\to H\ (n\to\alpha(n)\zeta)$ is such
that $R(\zeta)\lambda(n')=\alpha(n')R(\zeta)\ (\forall n'\in N$), so,
for any $\zeta,\eta\in H$, operator $R(\zeta)R(\eta)^*$ belongs to the
commutant $\alpha(N)'$ and $<\zeta,\eta>_N=R(\eta)^*R(\zeta)$
belongs to the commutant $\lambda(N)'$ which is isomorphic to the
opposite $C^*$-algebra $N^o$.

In a similar way, if $\beta$ be a faithful unital anti-$*$-representation 
of $N$ in $H$, then $H$ can be viewed simultaneously as a right $N$-module
$H_\beta$ or as a left $N^o$-module. An operator $R^o(\zeta):L_2(N)\to
H\ (n\to\beta(n^*)\zeta)$ is such
that $R^o(\zeta)\rho(n')=\beta(n'^*)R^o(\zeta)\ (\forall n'\in N$), so,
for any $\zeta,\eta\in H$, operator $R^o(\zeta)R^o(\eta)^*\in\beta(N)'$
and $<\zeta,\eta>_{N^o}=R^o(\eta)^*R^o(\zeta)\in\lambda(N)$.
\medskip

{\bf 1.1.2.} If $K$ is another Hilbert space carrying a faithful unital
$*$-repre\-sen\-ta\-tion $\gamma$ of $N$, we obtain a relative tensor
product $H\underset{N}{_\beta\otimes_\gamma} K$  from $H\otimes K$
introducing the
following inner product (and factorizing, if necessary):
$$
(\zeta_1\otimes\eta_1|\zeta_2\otimes\eta_2)=(\gamma(<\zeta_1,\zeta_2>_{
N^o})
\eta_1|\eta_2)\ (\forall \zeta_1,\zeta_2\in H, \eta_1,\eta_2\in K).
$$
Then there exists a linear operator $A_{\beta,\gamma}:H\otimes K\to
H\underset{N}{_\beta\otimes_\gamma} K$, the image of $\zeta\otimes\eta$
being denoted by $\zeta\underset{N}{_\beta\otimes_\gamma}\eta\
(\forall\zeta\in H,\eta\in K)$. One can also define
$H\underset{N}{_\beta\otimes_\gamma} K$ by introducing the following
inner product:
$$
(\zeta_1\otimes\eta_1|\zeta_2\otimes\eta_2)=
(\beta(<\eta_1,\eta_2>_{N}) \zeta_1|\zeta_2)\ (\forall
\zeta_1,\zeta_2\in H,
\eta_1,\eta_2\in K),
$$
this leads to a definition of a flip
$\sigma_{\beta,\gamma}:H\underset{N}{_\beta\otimes_\gamma} K\to
K\underset{N^o}{_\gamma\otimes_\beta} H$ by  $\sigma_{\beta,\gamma}:
\zeta\underset{N}{_\beta\otimes_\gamma}\eta\to\eta\underset{N^o}
{_\gamma\otimes_\beta}\zeta\ (\forall\zeta\in H,\eta\in K)$.
This allows to define a flip $\Sigma_{\beta,\gamma}:L(H\underset{N}
{_\beta\otimes_\gamma} K)\to L(K\underset{N^o}{_\gamma\otimes_\beta}
H)$ by
$\Sigma_{\beta,\gamma}(X)=\sigma_{\beta,\gamma}^*X\sigma_{\beta,\gamma}
\ (\forall X\in L(H\underset{N}{_\beta\otimes_\gamma} K))$.
Recall also that $\beta(n)\zeta\underset{N}{_\beta\otimes_\gamma}\eta=
\zeta\underset{N}{_\beta\otimes_\gamma}\gamma(n)\eta\
(\forall\zeta\in H,\eta\in K,n\in N)$.
\medskip

{\bf 1.1.3.} Finally, if $_\alpha H_\beta$ (resp., $_\gamma K_\delta$)
is an $(N,N)$-bimodule with respect to a $*$-representation $\alpha$
(resp., $\gamma$) and an anti-$*$-representation $\beta$ (resp., $\delta$),
then $H\underset{N}{_\beta\otimes_\gamma} K$ is also $(N,N)$-bimodule:
$n\cdot_\alpha[\zeta\underset{N}{_\beta\otimes_\gamma}\eta]\cdot_\delta
n'=[\alpha(n)\zeta]\underset{N}{_\beta\otimes_\gamma}[\delta(n')\eta]$
$(\forall\zeta\in H,\eta\in K,n,n'\in N)$.

The described construction does not depend on the choice of a faithful
trace on $N$ (up to an isomorphism of Hilbert spaces).

\end{subsection}

\begin{subsection}{Fibered product of $C^*$-algebras}
\label{fibered product}

Given two $C^*$-algebras $M_1,M_2$ such that $\beta(N)\subset
M_1,\gamma(N)
\subset M_2$, we construct their {\it fibered product} - a canonical
algebra
$(M_1)\underset{N}{_\beta{*}_\gamma} (M_2)$ isomorphic to the reduction
of $M_1\otimes M_2$ by some projection $e_{\beta,\gamma}$ depending on
$\beta,\gamma$.
\medskip

{\bf 1.2.1.} First, we define
$(M'_1)\underset{N}{_\beta\otimes_\gamma}(M'_2)$
as a $C^*$-algebra generated by $x\underset{N}{_\beta\otimes_\gamma} y$
$(\forall x
\in M'_1,y\in M'_2)$ and then $(M_1)\underset{N}{_\beta {*}_\gamma}
(M_2):= [(M'_1)\underset{N}{_\beta\otimes_\gamma}(M'_2)]'$. One can
verify that:
\begin{eqnarray*}
[(M_1)\underset{N}{_\beta{*}_\gamma} (M_2)]^o
&=&  ({M^o}_1) \underset{N^o}{_\beta{*}_\gamma} ({M^o}_2),\\
\Sigma_{\beta,\gamma}((M_1)\underset{N}{_\beta{*}_\gamma} (M_2))
&=& (M_2)\underset{N^o}{_\gamma{*}_\beta} (M_1),\\
(M_1\cap\beta(N)')\underset{N}{_\beta\otimes_\gamma}
(M_2\cap\gamma(N)') &\subset& (M_1)\underset{N}{_\beta{*}_\gamma}(M_2), \\
(M_1)\underset{N}{_\beta{*}_\gamma} (N)
= (M_1\cap\beta(N)')\underset{N}{_\beta\otimes_\gamma} 1_K
&=& (1_H)\underset{N}{_\beta\otimes_\gamma} (M_2\cap\gamma(N)').
\end{eqnarray*}

If $\beta'$ (resp., $\gamma'$) is an injective unital anti-$*$-homomorphism (resp., $*$-homomorphism) of $N$ into another $C^*$-algebra $L_1$ (resp., $L_2$), and $\Phi:M_1\to L_1$ (resp., $\Psi:M_2\to L_2$) is a unital $*$-homomorphism such that 
$\Phi\circ\beta=\beta'$ (resp., $\Psi\circ\gamma=\gamma'$), then
one can define a unital $*$-homomorphism
$\Phi\underset{N}{_\beta{*}_\gamma}\Psi:
(M_1)\underset{N}{_\beta{*}_\gamma} (M_2)\to
(L_1)\underset{N}{_{\beta'}\otimes_{\gamma'}} (L_2)$. This means, in particular, that if $\beta$ (resp.,$\gamma$) is an injective unital anti-$*$-homomorphism 
(resp., $*$-homomor\-phism) of $N$ into
$M_1$ (resp., $M_2$), then one can define
$(M_1)\underset{N}{_\beta{*}_\gamma}(M_2)$ (to be a $C^*$-algebra) without any references to a specific Hilbert space.

Obviously, if $\tau:P\to N$ is an isomorphism of $C^*$-algebras,
then one can identify
$H\underset{N}{_\beta\otimes_\gamma} K$ with
$H\underset{P}{_{\beta\circ\tau}\otimes_{\gamma\circ\tau}} K$,
$(M'_1)\underset{N}{_\beta\otimes_\gamma} (M'_2)$ with
$(M'_1)\underset{P}{_{\beta\circ\tau}\otimes_{\gamma\circ\tau}}
(M'_2),$ and $(M_1)\underset{N}{_\beta{*}_\gamma} (M_2)$ with
$(M_1)\underset{P}{_{\beta\circ\tau}{*}_{\gamma\circ\tau}} (M_2)$.
\medskip

{\bf 1.2.2.} Let now $H=K$ and $e_{\beta,\gamma}$ be the support of the 
operator $A_{\beta,\gamma}$ (see {\bf 1.1.2.}). Then one can show that
$e_{\beta,\gamma}\in\beta(N)\otimes\gamma(N)$ and
$e_{\beta,\gamma}(\beta(n)\otimes 1_H)=e_{\beta,\gamma}
(1_H\otimes\gamma(n))\ (\forall n\in N).$
Vice versa, let $M_1$ and $M_2$ be two $C^*$-algebras acting in $H$
such that $\beta(N)\subset M_1$ and $\gamma(N)\subset M_2$, and let
$e_{\beta,\gamma}\in\beta(N)\otimes\gamma(N)$ be a projection
satisfying the above relation. Then automatically $e
=(id\otimes\theta)(e_{\beta,\gamma})\in Z(\beta(N))$ and
$$||(e_{\beta,\gamma})^{1\over
2}(\zeta\otimes\eta)||^2=(\beta(<\eta,\eta>_N) e^{1\over
2}\zeta,e^{1\over 2}\zeta)\ (\forall \zeta,\eta\in H),$$
which equals to $||e^{1\over
2}\zeta\underset{N}{_\beta\otimes_\gamma}\eta||^2$. So $e^{1\over 2}$
is invertible and there exists an isometry
$V_{\beta,\gamma}:H\underset{N}{_\beta\otimes_\gamma} H\to H\otimes H$
with final support $e_{\beta,\gamma}$, defined by:
$V_{\beta,\gamma}(\zeta\underset{N}{_\beta\otimes_\gamma}\eta)=e_{\beta
,\gamma}
(e^{-{1\over 2}}\zeta\otimes\eta)\ (\forall \zeta,\eta\in H)$.

This isometry satisfies the following relations:
\begin{eqnarray*}
V_{\beta,\gamma}(x\underset{N}{_\beta\otimes_\gamma} y)
&=& (x\otimes y)V_{\beta,\gamma}\ (\forall
x\in\beta(N)',y\in\gamma(N)'); \\
(\beta(n)\otimes 1_H)V_{\beta,\gamma}
&=& (1_H\otimes\gamma(n)) V_{\beta,\gamma}\ (\forall n\in N); \\
(M'_1)\underset{N}{_\beta\otimes_\gamma} (M'_2)
&=& V^*_{\beta,\gamma}(M'_1\otimes M'_2)V_{\beta,\gamma}; \\
(M_1)\underset{N}{_\beta{*}_\gamma} (M_2)
&=& V^*_{\beta,\gamma}(M_1\otimes M_2)V_{\beta,\gamma},
\end{eqnarray*}
and $(M_1)\underset{N}{_\beta{*}_\gamma} (M_2)$ is isomorphic to the
reduced
$C^*$-algebra $(M_1\otimes M_2)_{e_{\beta,\gamma}}$. Let $i_{\beta,\gamma}$
be a non-unital injection of $(M_1)\underset{N}{_\beta{*}_\gamma} (M_2)$
into $M_1\otimes M_2$.
\end{subsection}

\begin{subsection}{Counital $C^*$-bialgebras}
\label{generalized Hopf}

\begin{definition}
\label{Def: generalized Hopf}
a) A pair  $(M,\Delta)$ is said to be a gen.\ $C^*$-bialgebra if $M$ is a
f-d $C^*$-algebra and a {\it coproduct} $\Delta:M\to M\otimes M$ is an
injective $*$-homomorphism
(not necessarily unital) such that  $(\Delta\otimes id)\circ\Delta
= (id\otimes\Delta)\circ\Delta$.

b) A triplet  ${\cal A}=(M,\Delta,S)$ is said to be a gen.
coinvolutive $C^*$-bialgebra if $(M,\Delta)$ is a gen.
$C^*$-bialgebra and a {\it co-involution} $S:M\to M$ is a unital
anti-$*$-automorphism such that $S^2=id$ and
$(S\otimes S)\circ\Delta=\varsigma\circ\Delta\circ S$,
where $\varsigma:M\otimes M\to M\otimes M$ is a usual flip.

c)  A homomorphism of gen.\ coinvolutive $C^*$-bialgebras
$(M_1,\Delta_1,$ $S_1)$ and $(M_2,\Delta_2,S_2)$ is a unital
$*$-homomorphism of $C^*$-algebras $\pi:M_1\to M_2$ such that
$(\pi\otimes\pi)\circ\Delta_1=\Delta_2\circ\pi,\ \pi\circ S_1
= S_2\circ\pi$.
Clearly, gen.\ coinvolutive $C^*$-bialgebras form a category.

d) \cite{Y} A quadruplet $GK= (M,\Delta,S,\phi)$ is said to be a
gen.\ Kac algebra if ${\cal A} = (M,\Delta,S)$ is a gen.
coinvolutive $C^*$-bialgebra and a {\it Haar trace} $\phi$ is a faithful
trace on $M$ such that $\phi\circ S=\phi$ and
$$
(id\otimes\phi)[(1\otimes y)\Delta(x)]
= S[(id\otimes\phi)(\Delta(y)(1\otimes x))]\ (\forall x,y\in M).
$$
\end{definition}

Gen.\ Kac algebra is {\it abelian}, i.e., $M$ is abelian (resp.,
{\it symmetric}, i.e., $\varsigma\circ\Delta=\Delta$), iff it is
isomorphic to the algebra of functions on  a finite groupoid (resp., to the
groupoid algebra of a finite groupoid) \cite{Y}.

Given a gen.\ $C^*$-bialgebra, $e=\Delta(1)$, let us define two sets:
$N_s=\{x\in M | \Delta(x)=e(1\otimes x)=(1\otimes x)e\}$ and
$N_t=\{x\in M |\Delta(x)=e(x\otimes 1)=(x\otimes 1)e\}$.

\begin{proposition}\cite{Val3} For any gen.\ $C^*$-bialgebra we have:
\label{Cartan subalgebras}
\begin{itemize}
\item[a)] $N_s$ and $N_t$ are $C^*$-subalgebras of $M$.
\item[b)] $N_s$ and $N_t$ commute: $nn'=n'n\ (\forall n\in N_s, n'\in
N_t)$.
\end{itemize}
\end{proposition}
{\it Proof}. a) It is clear that $N_s$ and $N_t$ are linear
$*$-invariant
subspaces of $M$ containing $1$. For any $x,y\in N_t$ we have:
\begin{eqnarray*}
\Delta(xy)
&=& \Delta(x)\Delta(y)=\Delta(x)e(y\otimes 1) \\
&=& \Delta(x)(y\otimes 1)= e(xy\otimes 1),
\end{eqnarray*}
so $N_t$ is a $C^*$-subalgebra of $M$; similarly for $N_s$.

b) For all $n\in N_s, n'\in N_t$ we have:
\begin{eqnarray*}
\Delta(nn')
&=& \Delta(n)e(n'\otimes 1)=\Delta(n)(n'\otimes 1) \\
&=& e(1\otimes n)(n'\otimes 1)=e(n'\otimes 1)(1\otimes n)\\
&=&\Delta(n')(1\otimes n)=\Delta(n'n),
\end{eqnarray*}
and the injectivity of $\Delta$ gives the result.

$N_s$ and $N_t$ are said to be {\em source and target Cartan
subalgebras}
respectively. In any generalized coinvolutive $C^*$-bialgebra 
$S(N_s)=N_t$.

\begin{definition}\label{counital algebra} A gen.\ coinvolutive
$C^*$-bialgebra is said to be {\it counital} if the {\it target
counital map}
$\eps_t:=\mu(id\otimes S )\Delta:M\to M$ (where $\mu:M\otimes M\to M$
is the multiplication in $M$) satisfies the following conditions:

1. $\eps_t(1)=1$.

2. $Im\eps_t\subset N_t$.

3. $(id\otimes\eps_t)\Delta(x)=e(x\otimes 1)\ (\forall x\in N_s)$.
\end{definition}
A definition of a morphism of generalized counital $C^*$-bialgebras is
obvious.
\begin{remark}\label{source counit}
a) Clearly, {\it source counital map} $\eps_s:=\mu(S\otimes
id)\Delta:M\to M$ has similar properties:

1'. $\eps_s(1)=1$.

2'. $Im\eps_s\subset N_s$.

3'. $(\eps_s\otimes id)\Delta(x)=(1\otimes x)e\ (\forall x\in N_t)$.

4. $S\circ\eps_t=\eps_s\circ S$.

b) From 3 and 3' we have $(id\otimes\eps_t)e=e=(\eps_s\otimes id)e$,
together
with 2 and 2' this gives $e\in N_s\otimes N_t$.

c) $(\Delta\otimes id)(e)=(id\otimes\Delta)(e)=(e\otimes 1)(1\otimes
e)=(1\otimes e)(e\otimes 1)$.

Indeed, let $e=\sum_{i=1}^n x_i\otimes y_i$ (where $x_i\in N_s, y_i\in
N_t$ - such a represen- tation is possible by b)). Then, using the
definition of $N_s$, we have:
$$(\Delta\otimes id)(e)=\sum_{i=1}^n \Delta(x_i)\otimes
y_i=\sum_{i=1}^n e(1\otimes x_i)\otimes y_i=(e\otimes 1)(1\otimes e).$$
The other equalities are now obvious.

d) One can easily verify that:
$$
\eps_t(x)^*=\eps_t(S(x)^*),\ \eps_t(xy)=\eps_t(x\eps_t(y)),\
\eps_t(xn)=\eps_t(xS(n))
$$
for all $x,y\in M,n\in N_t$, and similar relations for $\eps_s$.
\end{remark}

\begin{proposition}\label{bimodule map}
$\eps_t(nS(n')x)=n\cdot\eps_t(x)\cdot n'\ (\forall x\in M, n,n'\in
N_t).$
\end{proposition}
{\it Proof.} The definitions of $N_t, N_s$ and the properties of
$\Delta,\, S$ give:
\begin{eqnarray*}
\eps_t(nS(n')x)
&=& \mu(id\otimes S )\Delta(nS(n')x) \\
&=& \mu(id\otimes S )[(n\otimes 1)e(1\otimes S (n')) \Delta(x)] \\
&=& n\cdot\mu[(id\otimes S )\Delta(x)(1\otimes n')]=n\cdot\eps_t(x)\cdot n'.
\end{eqnarray*}

\begin{remark}
\label{source bimodule map}
a) Let us introduce the following structure of an $(N_t, N_t)$-bimodule
$_t M_s$ on $M$ : $n\cdot_t x\cdot_s n'=nS(n')x\ (\forall x\in M, n,
n'\in N_t)$.
Then Proposition \ref{bimodule map} means that $\eps_t: {}_tM_s\to N_t$
is the $(N_t, N_t)$-bimodule map.

b) One can prove the similar property for $\eps_s$:
$$\eps_s(xnS(n'))=S(n)\cdot\eps_s(x)\cdot S (n')\ (\forall x\in M,
n,n'\in N_t).$$

c) Proposition \ref{bimodule map} and axiom 1 imply that $\eps_t$
equals to identity on $N_t$ and to $S$ on $N_s$; similarly $\eps_s$
equals to identity on $N_s$ and to $S$ on $N_t$.
\end{remark}
\begin{proposition}\label{crucial} In any gen.\ counital
$C^*$-bialgebra
$$
e(S(n)\otimes 1)=e(1\otimes n),\ (S(n)\otimes 1)e=(1\otimes n)e\
(\forall n\in N_t).
$$
\end{proposition}
{\it Proof.} Clearly, the second relation can be deduced from the
first. Then, using successively the axiom 3, the definition of $N_t$,
Proposition \ref{bimodule map} and axiom 1, we have for any $ n\in
N_t$:
$$e(S(n)\otimes 1)=(id\otimes\varepsilon_t)\Delta(S(n))=
(id\otimes\varepsilon_t)e(1\otimes S(n))=e(1\otimes n).$$
\begin{definition}\label{Haar things} Let $(M,\Delta)$ be a gen.
counital $C^*$-bialgebra.

a) A projection $\Lambda\in M$ is said to be {\it a Haar projection} if
$$
x\Lambda=\eps_t(x)\Lambda\qquad(\forall x\in M), \quad S
(\Lambda)=\Lambda,\quad\eps_t(\Lambda)=1.
$$

b) A functional (resp., a faithful trace) $\phi$ on $M$ is said to be
{\it a Haar functional} (resp., {\it a Haar trace}) if
$$
(id\otimes\phi)\Delta=(\eps_t\otimes\phi)\Delta,\quad\phi\circ
S=\phi,$$
$\phi$ is said to be {\it normalized} if $(id\otimes\phi)e=1.$

c) A conditional expectation \cite{GHJ} 
$E_t : M\to N_t$ (resp., $E_s : M\to N_s$)
is said to be {\it target (resp., source) Haar conditional expectation}
on $(M,\Delta)$ if $(id\otimes E_t)\circ\Delta=\Delta\circ E_t$ (resp.,
$(E_s\otimes id)\circ\Delta=\Delta\circ E_s$).
\end{definition}
Obviously, the identities symmetric to a) and b) hold true for $\eps_s$.
\end{subsection}

\begin{subsection}{Counital Hopf bimodules}
\label{Hopf bimodules}

\begin{definition}\label{Def: Hopf bimodules} \cite{Val1}, \cite{EVal}.
A collection $B=(N,M,t,s,\Delta)$ is said to be a {\it Hopf bimodule},
if $N,M$ are  $C^*$-algebras, $s$ (resp., $t,\ \Delta$) is an injective
unital anti-$*$-homomorphism (resp., $*$-homomorphism) of $N\to M$
(resp., $N\to M,\  M\to M\underset{N}{_s{*}_t} M$) such that
$$
s(n)t(n')=t(n')s(n)\ \ \forall n,n'\in N,
$$
$$
\Delta(t(n))=t(n)\underset{N}{_s\otimes_t}(1),\qquad
\Delta(s(n))=(1)\underset{N}{_s\otimes_t} s(n)\ (\forall n\in N);
$$
$$
(\Delta\underset{N}{_s{*}_t} id)\circ\Delta
= ((id)\underset{N}{_s{*}_t}\Delta)\circ\Delta.
$$
\end{definition}

Let us clarify this definition. The structures of left, right
$N$-module and $(N,N)$-bimodule on $M$ (${}_t M, M_s$ and ${}_t M_s$
respectively), can be defined by:
$$
n\cdot_t m=t(n)m,\ m\cdot_s n'=ms(n'),\ n\cdot_t m\cdot_s n'=t(n)ms(n')
\ \ \forall m\in M, n,n'\in N.
$$
Then the fibered product $M\underset{N}{_s{*}_t} M$ is the $(N,N)$-bimodule
$$
{}_tM\underset{N}{{}_s{*}_t} M_s:\ n\cdot_t(m\underset{N}{_s{*}_t} m')\cdot_s
n'= (n\cdot_t m)\underset{N}{_s{*}_t}(m'\cdot_s n'),
$$
for all $m,m'\in M,\, n,n'\in N$.
So one can construct two fibered products:
$(M\underset{N}{_s{*}_t} M)\underset{N}{_s{*}_t} M$ and
$M\underset{N}{_s{*}_t}(M\underset{N}{_s{*}_t} M)$;
both of them are isomorphic to $p(M\otimes M\otimes M)p$, where
$p=(e_{s,t}\otimes 1)(1\otimes e_{s,t})
=(1\otimes e_{s,t})(e_{s,t}\otimes 1)$.
Two first conditions mean that $\Delta:_t M_s\to_t
M\underset{N}{_s{*}_t} M_s$
is an $(N,N)$-bimodule morphism; after that the last one is clear.

If $B=(N,M,t,s,\Delta)$ is a Hopf bimodule, then
$B^o=(N^o,M,s,t,\Sigma_{s,t}\circ\Delta)$ is also a Hopf bimodule; 
obviously,
$\Sigma_{s,t}\circ\Delta:M\to M\underset{N^o}{_t{*}_s}M$. If $N$ is
abelian and
$t=s$, then the two above Hopf bimodules coincide; we call such a Hopf
bimodule {\it symmetric}.

\begin{definition} \label{coinvolutive Hopf bimodules} Let
$B=(N,M,t,s,\Delta)$ be a Hopf bimodule. A {\it coinvolution} is such a
unital anti-$*$-automorphism  $S : M\to M$ that $S^2=id$, $S\circ t=s,$
and  $(S\underset{N}{_s{*}_t}S)\circ\Delta=\Sigma_{s,t}\circ\Delta\circ
S$, where $\Sigma_{s,t} : M\underset{N}{_s{*}_t}
M=M\underset{N^o}{_t{*}_s} M$
is a flip introduced in {\bf 1.1}, {\bf 1.2}. Then the collection
$B=(N,M,t,\Delta,S)$ is said to be a coinvolutive Hopf bimodule.
\end{definition}

Let us clarify this definition. The map 
$\Sigma_{s,t}\circ\Delta\circ S$ is a unital injective anti-$*$-homomorphism 
of $C^*$-algebras $M$ and $M\underset{N^o}{_t{*}_s} M$. Considering $S$ 
as an isomorphism between $M$ and $M^o$, the map 
$(S\underset{N}{_s{*}_t}S)$ is, according to {\bf 1.2}, the
isomorphism of $C^*$-algebras $M\underset{N}{_s{*}_t} M$ and
${M^o}\underset{N}{_t{*}_s} {M^o}$, the last algebra coincides with
$[M\underset{N^o}{_t{*}_s} M]^o$. Thus, the left-hand side of the
discussed equality is the unital injective anti-$*$-homomorphism of
$C^*$-algebras $M$ and $M\underset{N^o}
{_t{*}_s} M$.

\begin{definition} 
\label{Counital Hopf bimodules} 
A coinvolutive Hopf
bimodule $B=(N,M,t,\Delta,S)$ is said to be {\it counital} if the {\it
target counital map} $\epsilon_t:=\mu
[id\underset{N}{_s{*}_t}S]\circ\Delta : M\to M$ (where
$id\underset{N}{_s{*}_t}S : M\underset{N}{_s{*}_t}M\to M\underset{N}
{_s{*}_s}M^o$, the last one is a subalgebra in $M\otimes M^o$
(see {\bf 1.2}), and $\mu : M \otimes M^o\to M$ is the linear map given
by the multiplication in $M$) satisfies the following conditions:

1. $\epsilon_t(1)=1$.

2. $Im\, \epsilon_t\subset t(N)$.

3. $\epsilon_t: {}_tM_s\to t(N)$ is the $(N,N)$-bimodule map.

\end{definition}
A  definition of a morphism of counital Hopf bimodules is clear.

\begin{lemma}\label{unital equivalence} The categories of counital Hopf
bimodules and gen.\ counital $C^*$-bialgebras are equivalent.
\end{lemma}
{\it Proof.} a) Given a counital Hopf bimodule $B=(N,M,t,\Delta,S)$,
let us consider the projection $e_{s,t}\in s(N)\otimes t(N)$ such that
$e_{s,t}(S(n)\otimes 1)=e_{s,t}(1\otimes t(n))\ (\forall n\in N)$ and
$M\underset{N} {_s{*}_t} M= (M\otimes M)_{e_{s,t}}$, and the canonical
injection of $C^*$-algebras $i_{s,t} : M\underset{N} {_s{*}_t} M\to
M\otimes M$ (see {\bf 1.2}). Then $\tilde{\Delta} =i_{s,t}\circ\Delta :
M\to M\otimes M$ is an injective $*$-homomorphism (not necessarily
unital) such that:
$$(\tilde\Delta\otimes id)\circ\tilde\Delta=
(id\otimes\tilde\Delta)\circ\tilde\Delta,$$
$$\tilde\Delta(t(n))=(t(n)\otimes 1)\cdot e_{s,t},\
\tilde\Delta(s(n))=(1\otimes s(n))\cdot  e_{s,t},\ (\forall n\in N),$$
from where $\tilde\Delta(1)=e_{s,t}\in M\otimes M$. Now it is clear
that
$(M,\tilde{\Delta},S)$ is a gen.\ coinvolutive $C^*$-bialgebra
such that $\epsilon_t=\eps_t$ and $t(N)\subset N_t$, so the conditions
1 and 2 of Definition \ref{counital algebra} are satisfied. But the
definition of $N_t$, the structure of $\eps_t$ and the property
$\eps_t(1)=1$ imply that $\eps_t$ is identical on $N_t$; from
$Im\eps_t\subset t(N)$ we have $t(N)=N_t$, so $s(N)=N_s$.

Finally, using the above properties of $e_{s,t},\tilde\Delta$ and
$\epsilon_t$,
we have:
\begin{eqnarray*}
(id\otimes\eps_t)\tilde\Delta(s(n))
&=& (id\otimes\epsilon_t)e_{s,t}(1 \otimes s(n)) = e_{s,t}(1\otimes t(n)) \\
&=& e_{s,t}(s(n)\otimes 1)\ (\forall n\in N).
\end{eqnarray*}

b) Given a gen.\ counital $C^*$-bialgebra $(M,\tilde{\Delta},S)$
with Cartan
subalgebras $N_t$ and $N_s$, put $N=N_t, t=id_N, s=S$. Let $H=L_2(M)$
be the Hilbert space given by the GNS-construction for $(M,\phi)$,
where $\phi$ is a fixed faithful trace on $M$. Then, using the
properties of $e$ (Remark \ref{source counit} b), Proposition
\ref{crucial}), one can establish the following isomorphisms (see
\ref{relative tensor},\ref{fibered product}): $e(H\otimes H)\cong
H\underset{N}{_s\otimes_t} H$, $(e\otimes 1)(1
\otimes e)(H\otimes H\otimes H)\cong H\underset{N}{_s\otimes_t}
H\underset{N}{_s\otimes_t} H,\ (M'\otimes M')e\cong
M'\underset{N}{_s\otimes_t} M',\ (M'\otimes M'\otimes M')(e\otimes
1)(1\otimes e)\cong M'\underset{N}{_s\otimes_t}
M'\underset{N}{_s\otimes_t}  M'$ (spatially). Then
$M\underset{N}{_s{*}_t} M$ is isomorphic to the reduced $C^*$-algebra
$(M\otimes M)_e$ and $M\underset{N}{_s{*}_t}M\underset{N}{_s{*}_t}M$ is
isomorphic to $(M\otimes M\otimes M)_{(e\otimes 1)(1\otimes e)}$.

Taking into account the above isomorphisms, Definition \ref{counital
algebra} and Proposition \ref{bimodule map}, one can see that the
collection $(N,M,t,\Delta,S)$ is a counital Hopf bimodule. The
statements concerning morphisms are obvious.

Obviously, $e_{s,t}=1\otimes 1$ iff $N=\Bbb{C}$.

\begin{definition}
\label{Haar conditional expectation}
Let $B=(N,M,t,s,\Delta)$ be a  Hopf bimodule. A conditional expectation
$E_t:M\to t(N)$ (resp., $E_s:M\to s(N)$) is said to be a
{\it target (resp., source) Haar conditional expectation} on $B$ if:
$$
[(id)\underset{N}{_s{*}_t} E_t]\circ\Delta=\Delta\circ E_t, \quad
( {\text resp.,}
[(E_s)\underset{N}{_s{*}_t}(id)]\circ\Delta=\Delta\circ E_s).
$$
\end{definition}

Here $(id)\underset{N}{_s{*}_t} (E_t)$ is the restriction of the map
$id\otimes E_t:M\otimes M\to M\to t(N)$ to $M\underset{N}{_s{*}_t}M$;
since $E_t$ is a conditional expectation, this restriction is a well
defined map from $M\underset{N}{_s{*}_t}M=(M\otimes M)_{e_{s,t}}$ to
$M\underset{N}{_s{*}_t}t(N)=(M\otimes t(N))_{e_{s,t}}$. Similarly for
$(E_s)\underset{N}{_s{*}_t}(id)$.

Clearly $E$ is a source (resp., target) Haar conditional expectation on
a counital Hopf bimodule $B$ iff it is a source (resp., target) Haar
conditional expectation on the corresponding gen.\ counital $C^*$-bialgebra.
\begin{remark} \label{gen.Kac} There exists a duality theory for
gen.\ Kac algebras \cite{Y}: for any ${\cal K}=
(M,\Delta,S,\phi)$ one can construct a {\it dual} gen.\ Kac
algebra $\hat{\cal K}= (\hat
M,\hat\Delta,\hat S , \hat\phi)$ in such a way that the algebra dual
for $\hat{\cal K}$ is isomorphic to ${\cal K}$.
It is shown in \cite{Val3} that for any gen.\ Kac algebra:
a) $N_t=M\cap\hat M$ is a target Cartan subalgebra and the conditions
of Proposition 1.2 are satisfied; b) there exists the unique right Haar
conditional expectation $E : M\to N_t$, which can be defined as the
orthogonal projection
on $N_t$ in the Hilbert space given by the GNS-construction for
$(M,\phi)$ such that $\phi\circ E=\phi$.
\end{remark}
\end{subsection}
\end{section}


\begin{section} {Weak and generalized Kac algebras, Kac bimodules}

\begin{subsection}{Weak Kac algebras}\label{weak Kac algebras}

\begin{definition}
\label{Def: weak Kac algebra}
A {\it counit} for a gen.\ coinvolutive $C^*$-bialgebra $(M,$
$\Delta,S)$ is such a linear map $\eps: M\to \Bbb{C}$ that $(\eps\otimes id)
\Delta=(id\otimes\eps)\Delta=id$ and:
\begin{itemize}
\item[1).]
$\varepsilon(S(x))=\varepsilon(x),\ \varepsilon(x^*)
=\overline{\varepsilon(x)}.$
\item[2).]
$(\varepsilon\otimes\varepsilon)((x\otimes 1)e(1\otimes
y))=\varepsilon(xy)\ (x,y\in M,\ e:=\Delta(1)).$
\end{itemize}

A collection $WK=(M,\Delta,S,\varepsilon)$ satisfying 1)-2) and
\begin{itemize}
\item[3).]
$(\eps_s\otimes id)\Delta(x)=(1\otimes x)e,$
\end{itemize}
where $\eps_s=\mu(S\otimes id)\Delta,\ \mu:M\otimes M\to M$ is a multiplication
in $M$, is said to be a {\it weak Kac algebra.}
A homomorphism $\pi: WK_1\to WK_2$ of weak Kac algebras is a
homomorphism of their gen.\ coinvolutive $C^*$-bialgebras
(Definition~\ref{Def: generalized Hopf}~c)) which preserves the counits, i.e.,
$\eps_2\circ\pi= \eps_1$. Clearly, weak Kac algebras form a category.
\end{definition}
The above notion is the special case ($S^2=id$) of the notion of a
weak $C^*$-Hopf algebra introduced in \cite{BSz}.
It becomes a usual Kac algebra iff either $e=1\otimes 1$, or
$\varepsilon(xy)=\varepsilon(x)\varepsilon(y)$, or
$\mu(S\otimes id)\Delta(x)=\varepsilon(x)1\ (\forall x\in M)$.

\begin{proposition} \label{equivalence of axioms}
The set of axioms 2) and 3) is equivalent to the following set of
axioms:
\begin{itemize}
\item[A2.]
$(\varepsilon\otimes id)((x\otimes 1)e(1\otimes y))
=(\varepsilon\otimes id)((x\otimes 1)\Delta(y)),$
\item[A3.]
$(id\otimes\varepsilon\otimes id)[(e\otimes 1)(1\otimes\Delta(x))]
=e(1\otimes x),$
\item[A4.]
$\eps_s(x)=(id\otimes \varepsilon)((1\otimes x)e)\ (\forall x,y\in M).$
\end{itemize}
\end{proposition}
{\it Proof}. Indeed, applying $(id\otimes\varepsilon)$ to 3), we get A4.
This, in its turn, gives
\begin{itemize}
\item[A3'.]
$(id\otimes\varepsilon\otimes id)[(1\otimes\Delta(x)) (e\otimes 1)]
=(1\otimes x)e$,
\end{itemize}
which is clearly equivalent to A3.

Finally, using the above identities and the standard notation
$\Delta(x)=x_{(1)}\otimes x_{(2)}\ (\forall x\in M)$, we have 
equality A2':
\begin{eqnarray*}
(\varepsilon\otimes id)((1\otimes x)e(y\otimes 1))
&=& (\varepsilon\otimes\varepsilon\otimes id)
((1\otimes\Delta(x))(e\otimes 1)(y \otimes 1\otimes 1)) \\
&=& (\varepsilon\otimes\varepsilon)((1\otimes x_{(1)})e(y\otimes
1))x_{(2)} \\
&=& (\varepsilon\otimes id)(\Delta(x)(y\otimes 1)).
\end{eqnarray*}
which is equivalent to A2.
Vice versa, applying $\varepsilon$ to A2, we get 2), and we get 3)
combining A3' and A4.

The following relations are straightforward corollaries of A2-A4:
\begin{itemize}
\item[A3''.]
$(id\otimes\eps_t)\Delta(x)=e(x\otimes 1).$
\item[A4'.]
$\eps_t(x)=(\varepsilon\otimes id)(e(x\otimes 1))\ \ ({\text\ here}\
\eps_t=\mu(id\otimes S )\Delta).$
\end{itemize}

\begin{remark} 
\label{selfduality of axioms}
a) Let us equip the dual linear space $\hat M$ with the product and 
coproduct obtained by transposing the coproduct and product of $M$
by means of the canonical pairing $<,>:\hat M\times M\to\Bbb{C}$.
The unit of $\hat M$ is $\hat 1=\varepsilon$, the coinvolution $\hat S$
and the involution $*$ of $\hat M$ are defined by
$$
<\hat S (\alpha),x>=<\alpha,S(x)>,\quad
<\alpha^*,x>=\overline{<\alpha,S(x)^*>}\ (\forall x\in M,\alpha\in\hat M).
$$
One can show that $(\hat M,\hat\Delta, \hat S ,\hat\varepsilon)$ is also
a weak Kac algebra. Indeed, using linear maps $\eta: \lambda\to \lambda
e$
from $\Bbb{C}$ to $M\otimes M$ and $\varepsilon\circ\mu:M\otimes M\to
\Bbb{C}$, let us rewrite the above axioms as follows:
\begin{itemize}
\item[$\tilde{A2}$.]
$(\varepsilon\circ\mu\otimes\mu)(id\otimes\eta\otimes id)
= (\varepsilon\circ\mu\otimes id)(id\otimes\Delta).$
\item[$\tilde{A3}$.]
$(id\otimes\varepsilon\circ\mu\otimes id)(\eta\otimes\Delta)
= (id\otimes\mu)(\eta\otimes id).$
\item[$\tilde{A4}$.]
$\mu(S\otimes id)\Delta
= (id\otimes\varepsilon\circ\mu)(id\otimes\varsigma\otimes id)
(id\otimes\eta).$
\end{itemize}

Now it is clear that A2 and A3 are dual to one another, A4 is selfdual.

For the fact that $\hat M$ is a $C^*$-algebra see Corollary \ref{Dual Haar
trace}.

\noindent
b) The particular case of A3, in which $x=1$:
\begin{itemize}
\item[$A3^*$.]
$(id\otimes\varepsilon\otimes id)[(e\otimes 1)(1\otimes e)]=e$
\end{itemize}
is sufficient for getting A3. Indeed, using axioms A2, A3*, one has
$\forall x\in M$:
\begin{eqnarray*}
(id\otimes\varepsilon\otimes id)[(e\otimes 1)(1\otimes\Delta(x))]
&=& 1_{(1)}\otimes(\varepsilon\otimes id)[(1_{(2)}\otimes 1)\Delta(x)]
\\
&=& 1_{(1)}\otimes(\varepsilon\otimes id)[(1_{(2)}\otimes 1)e(1\otimes
x)] \\
&=& 1_{(1)}\otimes(\varepsilon\otimes id)[(1_{(2)}\otimes 1)e]\cdot x\\
&=& e(1\otimes x).
\end{eqnarray*}

\end{remark}

\begin{example}
As groups and their duals are trivial examples of usual Kac algebras,
groupoids and their duals give trivial examples of weak Kac algebras.
Let $G$ be a finite groupoid (see \cite{R} for definitions and terminology).
\newline
(a) $\Bbb{C}G$, the groupoid $C^*$-algebra of $G$, has a structure of a
cocommutative  weak Kac algebra given by
$$
\Delta(g)=g\otimes g,\quad S(g) =g^{-1}, \quad \eps(g)=1
\quad \text{ for all } g\in G.
$$
In this case $\eps_s(g) = s(g)=g^{-1}g$ and $\eps_t(g) =t(g)= gg^{-1}$ are
familiar source and target maps. The source and target
Cartan subalgebras coincide with $\Bbb{C}G^0$, where $G^0$
is the unit space of $G$.
\newline
(b) $\Bbb{C}(G)$, the $C^*$-algebra of complex-valued functions
on $G$, has a structure of a commutative  weak Kac algebra given by
$$
\Delta(\delta_g) =\sum_{xy=g} \delta_x\otimes\delta_y,\quad
S(\delta_g) = \delta_{g^{-1}},\quad \eps(\delta_g)=\delta_{g,t(g)}
= \delta_{g,s(g)} \quad \text{for all } g\in G.
$$
In this case $\eps_s(\delta_g) = \sum_{s(x)=g}\, \delta_x$ and
$\eps_t(\delta_g) = \sum_{t(x)=g}\, \delta_x$.  The source and
target Cartan subalgebras are
$N_s =\{f\in\Bbb{C}(G)\mid f(g)=f(s(g))\,\forall g\in G\}$
and  $N_t =\{f\in\Bbb{C}(G)\mid f(g)=f(t(g))\,\forall g\in G\}$.

The above weak Kac algebras are clearly dual to each other.
\end{example}

\begin{remark}
It follows from \cite{Y} and Theorem~\ref{weak2generalized}
below that every cocommutative (resp. commutative) weak Kac algebra
is isomorphic to $\Bbb{C}G$ (resp. $\Bbb{C}(G)$) for some groupoid $G$.
\end{remark}

Given two weak Kac algebras $M_1$ and $M_2$, one can construct their
tensor product $M_1\otimes M_2$ and direct sum $M_1\oplus M_2$ in an
obvious way.

Let us write $e=\sum_{i=1}^n x_i\otimes y_i$ with minimal possible $n$
(i.e., with both $\{x_i\}$ and $\{ y_i\}$ linearly independent).
Obviously,
one can choose $x_i$ and $y_i$ in such a way that $x_i=x^*_i,\
y_i=y^*_i.$

\begin{proposition}
\label{duality for Cartan subalgebras}
$\varepsilon(y_ix_k)=\varepsilon(x_ky_i)=\delta_{ik}$.
\end{proposition}
{\it Proof.} From $A3^*$ one gets:
$$
\sum_{i,j=1}^n x_i\otimes\varepsilon(y_ix_j)y_j = \sum_{i=1}^n x_i\otimes y_i,
\quad
\sum_{i,j=1}^n x_i\otimes\varepsilon(x_jy_i)y_j = \sum_{i=1}^n x_i\otimes y_i,
$$
so $\sum_{j=1}^n \varepsilon(y_kx_j)y_j =y_k =\sum_{j=1}^n
\varepsilon(x_jy_k)y_j,$
from where the result follows.
\medskip

Let us denote $N_s = \hbox{span}\{x_1,\dots x_n\}$ and
$N_t = \hbox{span}\{y_1,\dots y_n\}$. Then it is clear that
$e\in N_s\otimes N_t$, that both $N_s,N_t$ are $*$-invariant and
contain $1$.

\begin{proposition}
\label{properties of counits}
a) $\varepsilon_s$ (resp. $\varepsilon_t$) is a linear map from $M$ to
$N_s$ (resp. $N_t$); \quad
b) $\varepsilon_s\vert_{N_s} = id$, $\varepsilon_t\vert_{N_t} = id$.
\end{proposition}

{\em Proof.} a) is clear from A4 and A4'. b): $\varepsilon_s(x_k) =
\sum_{i=1}^n x_i\varepsilon(x_ky_i)
= x_k,$ similarly $\varepsilon_t(y_l) = y_l\ (\forall k,l=1,...,n)$.

\begin{proposition}
\label{Cartan properties}
$\Delta(\varepsilon_s(x))=(1\otimes\varepsilon_s(x))e, \quad
\Delta(\varepsilon_t(x)) =(\varepsilon_t(x)\otimes 1)e.$
\end{proposition}

{\it Proof.} The equality $(\Delta\otimes id)(e)=(id\otimes\Delta)(e)$
gives:
$$
\sum_{i=1}^n \Delta(x_i)\otimes y_i =\sum_{i=1}^n
x_i\otimes\Delta(y_i),
$$
from where, using Proposition~\ref{duality for Cartan subalgebras}
$$
\Delta(x_k)=\sum_{i=1}^n \Delta(x_i)\varepsilon(x_ky_i)
= \sum_{i=1}^n x_i\otimes(id\otimes\varepsilon)[(1\otimes
x_k)\Delta(y_i)]
\in N_s\otimes M.
$$
Similarly $\Delta(y_k)\in M\otimes  N_t$. Now Axioms 3), A3'' and
Proposition \ref{properties of counits} give the result.

\begin{corollary}
\label{target and source}
Any weak Kac algebra is a gen.\ counital $C^*$-bialgebra.
\end{corollary}

{\em Proof.} Relations 1 and 3 of Definition \ref{counital algebra}
follow from A4 and A4'. Since $N_s,N_t$ are $*$-invariant,
Proposition \ref{Cartan properties} shows that they belong to the
corresponding Cartan subalgebras, then Proposition \ref{properties of
counits} a) gives the relation 2 of Definition \ref{counital algebra}.

Moreover, $N_s$ and $N_t$ coincide with the corresponding Cartan
subalgebras, so they do not depend on the choice of a representation of
$e$ by means of $\{x_i,y_i\}_{i=1}^n$. Indeed, applying $\mu(S\otimes
id)$ to both sides of the equality in the definition of $N_s$, one can
see from A4 that $\eps_s$ is identical on this algebra. Similarly
$\eps_t$ is identical on $N_t$.

\begin{proposition}
\label{epsilon is a trace}
$\varepsilon(x)= Tr\,L_x\ (\forall x\in N_t)$,
where the linear operator $L_x$ acting in $N_t$ is defined by
$L_x(n)=xn\ (\forall n\in N_t),\ Tr\,L_x$ is its usual trace.
\end{proposition}
{\it Proof.} A4' shows that $\eps\circ\eps_t=\eps$, then from Corollary
\ref{target and source} and Remark \ref{source counit} d) we have for
all $x\in M,y\in N_t$: $\varepsilon(xy)=\varepsilon(yS(x))$.

Then in the representation
$(S\otimes id)e=\sum_{i=1}^n y'_i\otimes y_i=\sum_{i=1}^n
y_i\otimes y'_i$, where $y'_i=S(x_i),\ y_i\in N_t$, we
have from Proposition~\ref{duality for Cartan subalgebras} :
$\varepsilon(y'_k y_i)=\varepsilon(y_i y'_k)
= \delta_{i,k}\ (\forall i,k=1,...,n)$. This means that $\varepsilon$
is central on $N_t$ and that $\sum_{i=1}^n y'_i \cdot y_i=\sum_{i=1}^n
y_i\cdot y'_i=1$.
Now one can identify $N_t$ with its dual linear space by means of the
non-degenerated duality given by
$<y,z>=\varepsilon(yz)= \varepsilon(zy)\ (\forall y,z\in N_t)$.
Then, by the definition of a trace, for every $x\in N_t$ we have
(using the centrality of $\varepsilon$):
$Tr(L_x)=\sum_{i=1}^n <xy_i, y'_i>=\varepsilon(x(\sum_{i=1}^n y_i\cdot
y'_i))
=\varepsilon(x)$. In particular, $\eps(1)=\dim N_t$.

So the restriction of $\varepsilon$ to $N_t$ is a faithful trace.
On the other hand, using Axioms 1),2) of Definition~\ref{Def: weak Kac
algebra}, one can easily show that
$\varepsilon(\varepsilon_t(x)^*\varepsilon_t(y))=\varepsilon(x^*y)
\ (\forall x,y\in M)$, so $\varepsilon$ itself is a positive functional
on $M$.

\begin{proposition}
\label{formula for e}
Let $N_s \cong \oplus_{\alpha=1}^K\,M_{n_\alpha}({\Bbb{C}})$ and
$\{\,f_{pq}^\alpha\,\}_{p,q=1\dots n_\alpha}^{\alpha=1\dots K}$
be a system of matrix units in $N_s$. Then
$$
e=\sum\nolimits_\alpha\,\frac{1}{n_\alpha}\sum\nolimits_{pq}\,
f_{pq}^\alpha \otimes S(f_{qp}^\alpha).
$$
\end{proposition}

{\em Proof.}
Let us write  $e=\sum_{\alpha pq}\,f_{pq}^\alpha \otimes
g_{pq}^\alpha$, where $g_{pq}^\alpha \in N_t$. Then
$$
(f_{rr}^\beta \otimes 1) e (f_{ss}^\beta \otimes 1 ) = 
f_{rs}^\beta \otimes g_{rs}^\beta.
$$
On the other hand, using  Proposition~\ref{crucial} we have :
\begin{eqnarray*}
(f_{rr}^\beta \otimes 1) e (f_{ss}^\beta \otimes 1 ) 
&=& (1\otimes S(f_{rr}^\beta) ) e (1\otimes S(f_{ss}^\beta) )  \\
&=& \sum_{\alpha pq}\,f_{pq}^\alpha \otimes
S(f_{rr}^\beta) g_{pq}^\alpha S(f_{ss}^\beta), 
\end{eqnarray*}
from where $g_{rs}^\beta = c_{rs} S(f_{sr}^\beta)$
for some scalars $c_{rs}^\alpha$.
Proposition~\ref{duality for Cartan subalgebras} and
Remark~\ref{source counit}(d) imply that
$$
1= c_{rs}^\alpha \eps(f_{rs}^\alpha S(f_{sr}^\alpha)) =
c_{rs}^\alpha \eps(f_{rs}^\alpha f_{sr}^\alpha) =
c_{rs}^\alpha \eps(f_{rr}^\alpha),
$$
but $\eps(f_{rr}^\alpha) =n_\alpha$ by Proposition~\ref{epsilon is a
trace}, therefore $c_{rs}^\alpha = \frac{1}{n_\alpha}$, from where the
result follows.
\medskip

\begin{rem}
Morphisms of weak Kac algebras preserve Cartan subalgebras.
\end{rem}
Indeed, let $\rho : M \to \tilde{M}$ be a morphism and $\tilde{N_s}$,
$\tilde{N_t}$ be Cartan subalgebras of $\tilde{M}$.
Let us write $\Delta(1) = \sum_{i=1}^n\, x_i\otimes y_i$ with
$\{x_i\}_{i=1}^n$ linearly independent.
Then $N_s = \mbox{span}\{ x_i \}$ and since $\Delta(1_{\tilde{M}}) =
\sum_{i=1}^n\, \rho(x_i)\otimes \rho(y_i)$  we have $\tilde{N_s} =
\mbox{span}\{\rho(x_i)\}$, i.e., $\rho\vert_{N_s}$ is surjective.
On the other hand, by Proposition~\ref{epsilon is a trace}
$n=\dim N_s =\eps_M(1) = \eps_{\tilde{M}}(1) = \dim \tilde{N_s}$, so
$\rho\vert_{N_s}$ is injective and $N_s\cong \tilde{N_s}$.
Similarly, $N_t\cong \tilde{N_t}$.

\begin{corollary}
Weak Kac algebras with Cartan subalgebras isomorphic to the given
$C^*$-algebra $N$ form a subcategory of the category of weak Kac
algebras.
This subcategory is closed under duality.
\end{corollary}

We will show in Section 3 that any $C^*$-algebra $N$ can appear as
a Cartan subalgebra of some weak Kac algebra.

\end{subsection}

\begin{subsection}{The counital representation}\label{counital
representation}

\label{rep} Let $M = \oplus_{i\in\cal I}\,\Mdi $ be the $C^*$-algebra
of a weak Kac algebra, $P_i\ (i\in \cal I,\ \cal I$ is a finite set) be
the minimal central projections (i.e., selfadjoint idempotents) of $M$;
$\pi_i$ be the class of the irreducible representation $x\mapsto
P_ix\in\Mdi,\ x\in M$. Clearly, $\{\pi_i,\ i\in \calI\}$ is the set of
all classes of irreducible representations of $M$. Let $\chi_i$ be the
character of $\pi_i$.

Observe that the classes of non-degenerate representations of $M$ form
a ring $K_0(M)$ : if $\rho_i:M\to\calB(\calH_i),\ i=1,2$ are two
representations, then$$(\rho_1\oplus \rho_2)(x) := \rho_1(x)\oplus
\rho_2(x)
\quad\mbox{and}\quad(\rho_1\times\rho_2)(x) :=
(\rho_1\otimes\rho_2)\Delta(x)$$
are representations of $M$ in Hilbert spaces $\calH_1\oplus\calH_2$ and
$(\rho_1\otimes\rho_2)e(\calH_1\otimes\calH_1)$ respectively.
The set of all irreducible classes $\{\pi_i,\ i\in \calI\}$ forms a
linear basis of  $K_0(M)$. $K_0$-ring for usual Kac algebras was studied in
\cite{N}.

The class of a representation $\rho$ is completely determined by its
(normalized) character $\chi_\rho$, the product and the direct sum of
representations correspond to the  product and the direct sum of
characters
respectively. Let $\rho^*$ be the representation corresponding to the
character $\chi_\rho\circ S$, then the map $\rho\mapsto\rho^*$
defines an antimultiplicative involution in $K_0(M)$.

Consider the {\it counital representation} $\pi_\eps$ of $M$ associated
by
the GNS-con- struction with the positive functional $\eps$.
Proposition~\ref{epsilon is a trace} shows that $\eps$ is faithful on
$N_t$
and we also have
$\eps((\eps_t(x)-x)^*(\eps_t(x)-x))=0\ (\forall x\in M)$, therefore,
$\pi_\eps$ acts on the Hilbert space $N_t$ equipped with a scalar
product
$(x,\ y):=\eps(y^*x)$ in the following way:
$$
\pi_\eps(x)\eps_t(y) = \eps_t(xy) \quad\mbox{ for all}\quad x,y\in M.
$$

\begin{proposition}(cf. \cite{BSz})
\label{unit}
The class of $\pi_\eps$ is a unit for $K_0(M)$, i.e., $\pi_\eps\times
\rho$ and $\rho\times\pi_\eps$ are equivalent to $\rho$
for any representation $\rho$.
\end{proposition}

{\em Proof.}
Let us compute the character $\chi_\eps =\Tr\pi_\eps$.
Let $e=\sum_k x_k\otimes y_k$, where the set $\{y_k\}$ forms a basis of
$N_t$. Therefore, for all $x\in M$ we have
$$
\pi_\eps(x)y_k = \eps_t(xy_k) = \sum_i\,\eps(x_ixy_k)y_i,
$$
and
\begin{eqnarray*}
\chi_\eps(x)
&= & \sum\nolimits_k\,\eps(x_kxy_k)=\eps(\mu(e(x\otimes 1))) \\
&= & (\eps\otimes\eps)((1\otimes x\1)e(x\2\otimes 1))=\eps(\mu\circ
\Delta(x)).
\end{eqnarray*}
Using this formula we compute
\begin{eqnarray*}
(\chi_\eps\otimes\chi)\Delta(x)
&= & (\eps\otimes\chi)((\mu\otimes id)(\Delta\otimes id)\Delta(x))\\
&= & (\eps\otimes\chi)((x\1\otimes 1)e(1\otimes x\2))=\chi(x),
\end{eqnarray*}
for all characters $\chi$. Similarly, one can show that
$(\chi\otimes\chi_\eps)\Delta(x) = \chi(x)$.
Hence, $\pi_\eps\times \rho$, $\rho\times\pi_\eps$, and $\rho$ are
equivalent for all representations $\rho$.
\medskip

As a representation of a f-d $C^*$-algebra, $\pi_\eps$ is equivalent to
the sum of irreducible representations with some multiplicities:
$$
\pi_\eps\simeq\oplus_{i\in\calS}\,\nu_i\pi_i,
$$
where $\calS\subset\calI$ and $\nu_i\geq1$ for all $i\in \calS$.

\begin{proposition}(cf. \cite{BSz}, 2.4.)
\label{mult}
The representation $\pi_\eps$ is multiplicity free, i.e.
$\nu_i=1$ for all $i\in\calS$. We have $\pi_i^*= \pi_i$
for any $i\in\calS$ (and so $\pi_\eps^* = \pi_\eps$).
\end{proposition}

{\em Proof.}
Since $\pi_\eps$ is the unit of $K_0(M)$, we have $\pi_\eps^* =
\pi_\eps$
and $\pi_j = \pi_j\times\pi_\eps
=\sum_{i\in\calS}\,\pi_j\times\nu_i\pi_i$
for all $j\in\calI$. In particular, the right-hand side of the last
equality
must be irreducible. This is the case iff there exists a unique index
$u(j)\in \calS$ such that $\nu_{u(j)}=1$ and $\pi_j\times\pi_i
= \delta_{i\,u(j)}\, \pi_j$ for all $i\in\calS$.

On the other hand, for any irreducible $\pi_j$ one has
$\pi_j\times\pi_j^*\neq 0$. Indeed, if $(P_j\otimes S (P_j))e=0$,
then  applying $\mu(id\otimes S)$ we get $P_j=0$, a contradiction.
This implies $\pi_{u(j)}=\pi_j^*$ and, hence, $\pi_j\times\pi_j^*
=\pi_j$ for all $j\in\calS$. Applying the involution to the last
equality
we get $\pi_j = \pi_j^*$ and $\nu_{u(j)}= \nu_j =1\ (\forall
j\in\calS)$.

\begin{remark}\label{support}
Let $p_\eps$ be the support of $\eps$, i.e., a non-zero projection
$p$ minimal with respect to the property $\eps(xp)=\eps(px)=\eps(x)\
(\forall x\in M)$ (which is equivalent to
$\eps_s(px)=\eps_s(x)$ and $\eps_t(xp)=\eps_t(x)$).
Proposition~\ref{mult} shows that $p_\eps=\sum_{i\in\calS}\,p_i$,
where $p_i$ are minimal in $M$ projections whose central supports $P_i$
are mutually orthogonal. Obviously $\eps\circ S=\eps$ implies
$S(p_\eps)=p_\eps$.

Since $\pi_\eps$ is multiplicity free, we have
$$\dim N_t = \dim\calH_\eps =
\sum\nolimits_{i\in\calS}\,(\dim P_iM)^{\frac{1}{2}} =
\sum\nolimits_{i\in\calS}\,\dim Mp_i = \dim Mp_\eps.$$
Therefore, the map $\varepsilon_t: Mp_\eps\to N_t$ is a linear
isomorphism of
vector spaces. In other words, $\{x\in M| \eps_t(x)=0\}
= \{x\in M|xp_\eps=0\}$, this implies $(x-\eps_t(x))p_\eps = 0$ and
$xp_\eps =\eps_t(x)p_\eps\ (\forall x\in M)$.
Replacing $x$ by $S(x)$ and applying $S$ to the last equality,
we also get $p_\eps x= p_\eps\eps_s(x)$.
\end{remark}

\begin{proposition}
\label{Ideals}
Let $I_s =\{\,y\in M\mid yx = y\eps_s(x),\, \forall x\in M\,\}$
and $I_t =\{\,y\in M\mid xy = \eps_t(x)y,\, \forall x\in M\,\}$.
Then $I_s\cap I_t = \{\,y\in M \mid y = p_\eps y p_\eps\,\}$.
\end{proposition}

{\em Proof.}
Clearly, $I_s$ is a left ideal in $M$, so there exists a projection
$p\in M$
such that $I_s = Mp$. Since $p_\eps\in I_s$, we have
$p_\eps= p_\eps p = p_\eps p p_\eps$ and $p' = p - p_\eps \in I_s$ is a
projection orthogonal to $p_\eps$. Hence,
$$
p' = {p'}^2 = p'\eps_s(p') = p'\eps_s(p_\eps p') = 0.
$$
So $ p= p_\eps$ and $I_s = Mp_\eps$, similarly $I_t = p_\eps M$; these
equalities give the result.

\begin{theorem}
\label{Haar Theorem} 
Given a weak Kac algebra, there exists a unique Haar projection
$\Lambda=p_\eps\in M$.
\end{theorem}

{\em Proof.}
Proposition~\ref{Ideals} implies that any $\Lambda$ satisfying
\ref{Haar things} a) has the form $\Lambda= p_\eps \Lambda p_\eps =
\sum_{i\in\calS}\,\lambda_i p_i$ for some scalars $\lambda_i$; then,
due to
the above relations,
$\eps(p_i) =\eps(\Lambda p_i) = \lambda_i\eps(p_i)$ from where
$\lambda_i =1$
and $\Lambda = p_\eps$. Since $p_\eps$ satisfies \ref{Haar things} a),
the proof is completed.

Now let us introduce and study a {\it counital quotient} of a weak Kac
algebra arising from its counital representation $\pi_\eps$.

\begin{proposition}
\label{Canonical quotient}
\begin{enumerate}
\item[(a)]
Let $\{\pi_i\}_{i\in\calS}$ be the set of irreducible representations
of $M$ contained in the decomposition of $\pi_\eps$. Then for any
$i\in\calS$ $M_i=P_iM$ is a weak Kac algebra with a comultiplication
$\Delta_i(x)= (P_i\otimes P_i)\Delta(x)$, an antipode $ S_i(x) =
 S(x)$, and a counit $\varepsilon_i(x)= \varepsilon(x)$.
\item[(b)]
Let $P_\eps=\sum_{i\in\calS} P_i$ be the central support of $\eps$.
Then $x \mapsto P_\eps x$ is a surjective morphism of weak Kac algebras
$M$ and $M_\eps =\oplus_{i\in\calS}\,M_i$, i.e., $M_\eps$ is a quotient
of $M$.
\end{enumerate}
\end{proposition}

{\em Proof.}
Proposition~\ref{mult} implies that $\pi_i = \pi_i^*$ and
$\pi_i\times\pi_j =\delta_{ij}\pi_i$ for all $i,\,j\in\calS$.
This means that $(P_i\otimes P_i)\Delta = (P_i\otimes P_i)\Delta\circ P_i$
and  $ S\circ P_i=P_i\circ S$. In particular,
$\Delta_i(x) \in P_iM\otimes P_iM$ and $ S_i(x) \in P_iM\ (\forall x\in
P_iM)$. We also have that $(P_\eps\otimes P_\eps)\Delta =\Delta_i$.
$x\in P_iM$.

Clearly, $(P_i M,\,\Delta_i,\, S_i)$ is a gen.\ coinvolutive
$C^*$-bialgebra. Due to the definitions of $\Delta_i, S_i, \eps_i$
and the centrality of $P_i$, $\eps_i$ is a counit for $P_i M$
satisfying the axioms of Definition \ref{Def: weak Kac algebra},
which proves (a). Obviously, $M_\eps =\oplus_i\, M_i$ is
a weak Kac algebra. To see that $x\mapsto P_\eps x$ is a (clearly
surjective) morphism of weak Kac algebras, it suffices to note that
$(P_\eps \otimes P_\eps)\Delta(x) =\sum_i\Delta_i(P_\eps x)$,
$P_\eps S(x) = \sum_i S_i(P_\eps x)$, and
$\eps(x) = \sum_i\eps_i(P_\eps x)\ \forall x\in M$.
\medskip

\begin{define}
The weak Kac algebra $M_\eps$ defined in Proposition~\ref{Canonical
quotient} is said to be a counital quotient of $M$.
\end{define}

By Proposition~\ref{Canonical quotient}(a), $M_\eps$ is a direct sum
of weak Kac algebras which are simple as algebras. We will
classify all such algebras in Section 3.
\medskip

\begin{rem}
$M_\eps$ is a minimal quotient weak Kac algebra of $M$, i.e. if
$\rho : M \to \tilde{M}$ is a surjective morphism of weak Kac algebras,
then
$\rho(M_\eps)\cong \tilde{M_\eps}$.
\end{rem}
Indeed, $\eps(x) =\tilde{\eps}(\rho(x))\ (\forall x\in M)$ implies
$\rho(P_\eps)=\tilde{P_\eps}$, so $\rho\vert_{M_\eps} : M_\eps\to
\tilde{M_\eps}$ is surjective.
If it is not injective then $\rho(P_i)=0$ for some $i\in\calS$, which
contradicts to $\eps(P_i)\neq 0$. Therefore, $\rho\vert_{M_\eps}$ is an
isomorphism.

\end{subsection}


\begin{subsection}{Haar traces}

Let us consider representations $L_x: y\to xy,\ R^*_\alpha: y\to(id\otimes\alpha)
\Delta(y)\ (\forall x,y\in M,\alpha\in\hat M)$ \cite{M} of the algebras $M$
and $\hat M$ respectively in $End(M)$.

\begin{proposition}
\label{RL}
a) $R^*_\alpha L_x= \alpha_{(1)}(x_{(2)})L_{x_{(1)}}R^*_{\alpha_{(2)}},
\quad
b)\ R^*_{\hat\varepsilon_t(\alpha)}=L_{(id\otimes\alpha)e},$
where
$\hat\varepsilon_t(\alpha) = (\hat\varepsilon\otimes id_{\hat M})
(\hat\Delta(\varepsilon)(\alpha\otimes\varepsilon)).$
\end{proposition}

{\em Proof.}
a) We have
\begin{eqnarray*}
R^*_\alpha L_x y
&=& (id\otimes\alpha)\Delta(xy) = x_{(1)}y_{(1)}\alpha(x_{(2)}y_{(2)}) \\
&=&  x_{(1)}y_{(1)}\hat\Delta(\alpha)(x_{(2)}\otimes y_{(2)})
=x_{(1)}y_{(1)}\alpha_{(1)}(x_{(2)})\alpha_{(2)}(y_{(2)}) \\
&=& \alpha_{(1)}(x_{(2)}) L_{x_{(1)}}R^*_{\alpha_{(2)}}y.
\end{eqnarray*}

b) Clearly the map $\hat\varepsilon_t$ is dual to the map
$\varepsilon_t$, and we have, using A3'':
\begin{eqnarray*}
R^*_{\hat\varepsilon_t(\alpha)}(x)
&=& (id\otimes\hat\varepsilon_t(\alpha))\Delta(x)
 = (id\otimes\alpha\circ\varepsilon_t)(\Delta(x))\\
&=& (id\otimes\alpha)(e(x\otimes 1)) = L_{(id\otimes\alpha)e}(x).
\end{eqnarray*}

\begin{corollary}
\label{RL-corollary}
$R^*_{\alpha_{(1)}} L_x R^*_{\hat S (\alpha_{(2)})}=
L_{x_{(1)}}\alpha(x_{(2)}).$
\end{corollary}

{\em Proof.} Indeed, using the notation
$(\hat\Delta\otimes id_{\hat M})\hat\Delta(\alpha)
=(id_{\hat M}\otimes\hat\Delta)\hat\Delta(\alpha)=
\alpha_{(1)}\otimes\alpha_{(2)}\otimes\alpha_{(3)}$
and Proposition~\ref{RL}, we have:
\begin{eqnarray*}
R^*_{\alpha_{(1)}} L_x R^*_{\hat S (\alpha_{(2)})}
&=& L_{x_{(1)}}\alpha_{(1)}(x_{(2)})R^*_{\alpha_{(2)}\hat S
(\alpha_{(3)})}
=L_{x_{(1)}}\alpha_{(1)}(x_{(2)})R^*_{\hat\varepsilon_t(\alpha_{(2)})}
\\
&=& L_{x_{(1)}}\alpha_{(1)}(x_{(2)})L_{(id\otimes\alpha_{(2)})e}=
L_{x_{(1)}}\alpha_{(1)}(x_{(2)})L_{1_{(1)}}\alpha_{(2)}(1_{(2)}) \\
&=& L_{x_{(1)}1_{(1)}}\alpha(x_{(2)}1_{(2)})
= L_{x_{(1)}}\alpha(x_{(2)}).
\end{eqnarray*}

\begin{proposition}
Let us consider the following faithful trace on $M$:
$$\theta(x)=Tr(L_x),$$
where $Tr$ is the usual trace of a linear operator. Then we have:
$$
(\theta\otimes id)(\Delta(x))
= (\theta\otimes S )(e(x\otimes 1))
=(\theta\otimes\varepsilon_s\circ S)(\Delta(x))\ (\forall x\in M).
$$
\end{proposition}

{\em Proof.}
Using Proposition~\ref{RL}, Corollary~\ref{RL-corollary}
and the properties of $Tr$, one has for any $x\in M,\alpha\in\hat M$:
\begin{eqnarray*}
(\theta\otimes\alpha)(\Delta(x))
&=& Tr[L_{x_{(1)}}\alpha(x_{(2)})]
  = Tr[R^*_{\alpha_{(1)}} L_xR^*_{\hat S (\alpha_{(2)})}] \\
&=& Tr[R^*_{\hat S (\alpha_{(2)})\alpha_{(1)}} L_x ]
  = Tr[R^*_{\hat\varepsilon_t(\hat S (\alpha))}L_x] \\
&=& Tr[L_{(id\otimes\hat S (\alpha))e}L_x]
  = Tr[L_{(id\otimes\alpha\circ S) (e(x\otimes 1))}]\\
&=& (\theta\otimes\alpha\circ S)(e(x\otimes 1)).
\end{eqnarray*}

\begin{remark}
\label{H}
Since $(\theta\otimes id)(\Delta(x))\in N_s$, $\theta$ is a
Haar trace ($\theta\circ S=\theta$ from $S^2=id$ - see \cite{KP}, Lemma
3.1), but generally it is not normalized.
\end{remark}

\begin{theorem}
\label{weak2generalized}
Let $(M,\Delta,S,\varepsilon)$ be a weak Kac algebra and $\phi$ a Haar trace on
it. Then the collection $(M,\Delta,S,\phi)$ is a gen.\ Kac algebra.
\end{theorem}

{\em Proof.} It is enough to show that
$$
(\phi\otimes id)[(y\otimes 1)\Delta(x)]
 = S[(\phi\otimes id)(\Delta(y)(x\otimes 1))]\ (\forall x,y\in M).
$$
Indeed, using successively A3'', A4' and Definition of $\eps_t$ one has:
\begin{eqnarray*}
(\phi\otimes id)[(y\otimes 1)\Delta(x)]
&=& (\phi\otimes id) [\Delta(x)e(y\otimes 1)] \\
&=& (\phi\otimes id)[\Delta(x) (id\otimes\varepsilon_t)\Delta(y)]\\
&=& (\phi\otimes id)
[\Delta(x)(id\otimes\mu(id\otimes S ) \Delta)\Delta(y)] \\
&=& (\phi\otimes id)[\Delta(xy_{(1)})]S(y_{(2)}) 
\end{eqnarray*}

Then we compute, for any $z\in M$, using the definition of a Haar trace, relation
$S\circ\eps_t\circ\eps_s = \eps_s$ and A3'':
\begin{eqnarray*}
(\phi\otimes id)\Delta(z) = (\phi\otimes\eps_s)\Delta(z)
&=& (\phi\otimes[S\circ\eps_t\circ\eps_s])\Delta(z)\\
&=& [S\circ\eps_t](\phi\otimes id)\Delta(z)
= S[(\phi\otimes\eps_t)\Delta(z)]\\
&=& S[(\phi\otimes\eps_t)(e(z\otimes 1)]
\end{eqnarray*}

Finally, using this relation, we have:
\begin{eqnarray*}
(\phi\otimes id)[\Delta(xy_{(1)})]S(y_{(2)})
&=& S(\phi\otimes id)[e(xy_{(1)}\otimes 1)]S(y_{(2)}) \\
&=& S[(\phi\otimes id)(\Delta(y)(x\otimes 1))].
\end{eqnarray*}

Now let us compute $\Delta(p_\varepsilon)$.

\begin{lemma}
\label{Rank 1}
Let $A=M_d(\Bbb{C})$, $B=M_d(\Bbb{C})$, and $S:A\to B$
be an antiautomorphism of algebras. Let $Q\in A\otimes B$ be a
projection
of rank 1 in $A\otimes B$ such that $\mu(S\otimes\mbox{id}{}_B)Q = 1_B$
and $\mu(\mbox{id}{}_A\otimes S ^{-1})Q = 1_A$.
Then $Q=\frac{1}{d}\sum_{km}\,e_{km}\otimes S (e_{mk})$ for any system
$\{\,e_{km}\,\}_{k,m=1}^d$ of matrix units in $A$.
\end{lemma}

{\em Proof.} If $\,\{e_{km}\,\}_{k,m=1}^d$ is a system of matrix units
in $A$, then $\{\,S(e_{nl})\,\}_{l,n=1}^d$ is a system of matrix units
in $B$.
Since $Q$ has rank 1, we can write it as
$$Q = \sum_{klmn}\,\gamma_{kl}\beta_{mn}\,e_{km}\otimes S (e_{nl}),$$
for some $\beta_{mn},\ \gamma_{kl}\in\Bbb{C}\ (k,l,m,n =
1\dots d)$. By the hypothesis,
$$
1_B =  \mu(S\otimes\mbox{id}{}_B)Q =  \sum_{klmn}\,\gamma_{kl}\beta_{mn}\,
S(e_{nl}e_{km})= (\sum_{k}\,\gamma_{kk})\,\sum_{mn}\,\beta_{mn}\,S(e_{nm}),
$$
so $\beta_{mn}=\delta_{mn}\,\beta$ for some $\beta\in\Bbb{C}$.
Similarly, $\gamma_{kl}=\delta_{kl}\,\gamma$ for some
$\gamma\in\Bbb{C}$.
Thus, $ Q=\gamma\beta\,\sum_{km}\,e_{km}\otimes S (e_{mk})$
and condition $Q^2=Q$ gives $\gamma\beta=\frac{1}{d}$.

\begin{lemma}
\label{Preliminary}
Let $P_\eps = \sum_{i\in\calS}\,P_i$ be the central support
of $p_\eps$. Then
$$
(P_\eps\otimes P_\eps)\Delta(p_\eps)=
\sum_{i\in\calS}\,\frac{1}{d_i}\sum_{kl}\,
e_{kl}^{(i)}\otimes S (e_{lk}^{(i)}),
$$
where $\{\,e_{kl}^{(i)}\,\}{}_{k,l=1}^{d_i}$ is any system of matrix
units
in $P_iM$.
\end{lemma}

{\em Proof.}
The proof of Proposition~\ref{mult} shows that $P_i=S(P_i),\,i\in\calS$
and $ (P_\eps\otimes P_\eps)\Delta(p_\eps) = \sum_{i\in\calS}\,Q_i, $
where $Q_i = (P_i\otimes P_i)\Delta(p_\eps)$ is a projection
of rank 1 in  $P_iM\otimes P_iM$ such that $\mu(S\otimes id)Q_i
=P_i\eps_s(p_\eps) = P_i$. Similarly $\mu(id\otimes S )Q_i=P_i$.
Applying Lemma~\ref{Rank 1} we get
$ Q_i = \frac{1}{d_i}\sum_{kl}\,e_{kl}^{(i)}\otimes S (e_{lk}^{(i)}),$
from where the result follows.

\begin{lemma}
\label{theta}
$(\theta\otimes id)\Delta(p_\eps) = (id\otimes\theta)\Delta(p_\eps) =
1$.
\end{lemma}

{\em Proof.}
Let us choose a system of matrix units
$\{\,e_{kl}^{(i)}\,\}{}_{k,l=1}^{d_i}$
in $P_\eps M$ in such a way that $p_i = e_{11}^{(i)}$ and
$S(e_{lk}^{(i)}) = e_{kl}^{(i)}$ for all $i,\,k,\,l$.
We have, using Remark~\ref{H}, relation
$\Delta(p_\eps)=(P_\eps\otimes P_\eps)\Delta(p_\eps)$ and
Lemma~\ref{Preliminary}:
\begin{eqnarray*}
(\theta\otimes id)\Delta(p_\eps)
&=& (\theta\otimes\eps_s)\Delta(p_\eps)
= (\theta\otimes\eps_s)((1\otimes p_\eps)\Delta(p_\eps)) \\
&=& \sum_{i\in\calS}\,\frac{1}{d_i}\,
     \sum_{kl}\,\theta(e_{lk}^{(i)})\,\eps_s(e_{11}^{(i)}e_{kl}^{(i)})
=  \sum_{i\in\calS}\,\frac{1}{d_i}\,
     \sum_{kl}\, d_i\,\delta_{kl}\,\eps_s(\delta_{1k}e_{1l}^{(i)}) \\
&= & \sum_{i\in\calS}\,\eps_s(e_{11}^{(i)})=\eps_s(p_\eps)=\eps_s(1) =1.
\end{eqnarray*}
The second part is similar.

Let us generalize now the formula known for Kac algebras (\cite{ES},
6.3.7).

\begin{proposition}
\label{Evaluation}
For any system $\{\,e_{kl}^{(i)};\; k,l=1, \dots, d_i,\,i\in\calI\,\}$
of matrix units in $M$
$$
\Delta(p_\eps) =
\sum_{i\in\calI}\,\frac{1}{d_i}\sum_{kl}\,
e_{kl}^{(i)}\otimes S (e_{lk}^{(i)}).
$$
In particular, $\varsigma\Delta(p_\eps)=\Delta(p_\eps)$.
\end{proposition}

{\em Proof.} Let us write $\Delta(p_\eps) = \sum_{ij\in\calI}\,R_{ij}$,
where $R_{ij}$ is a projection in $P_iM\otimes P_jM$. From Lemma~\ref{theta} we
have
\begin{eqnarray*}
1 &=& (id\otimes\theta)\Delta(p_\eps) =
   \Sigma_i\,\left(\Sigma_j\,d_j(id\otimes\chi_j)R_{ij})\right) P_i \\
1 &=& (\theta\otimes id)\Delta(p_\eps)
   = \Sigma_j\,\left(\Sigma_i\,d_i(\chi_i\otimes id)R_{ij})\right) P_j.
\end{eqnarray*}
Applying $\chi_i$ to the first equality and $\chi_j$ to the second one,
we get :
$$
d_i = \Sigma_j\, d_j\,\rank(R_{ij}), \qquad
d_j = \Sigma_i\, d_i\,\rank(R_{ij}).
$$
Note that $\eps_s(p_\eps)=1$ implies that $R_{ii^*}\neq 0$
(we can assume that the involution acts on the set $\calI$).
Using this fact and the above relations, we conclude that
$\rank(R_{ij})=\delta_{ij^*}$. Thus, $\mu(S\otimes id)R_{ii^*}=
\mu(id\otimes S )R_{ii^*}=P_i$
and application of Lemma~\ref{Rank 1} completes the proof.

\begin{corollary}
\label{Dual Haar trace}
$p_\eps=\hat\phi_\eps$ is a normalized Haar trace on $\hat M$. $(\hat
M,\ \hat\Delta,\ \hat S,\ \hat\eps)$ is a weak Kac algebra.
\end{corollary}

{\em Proof.} Clearly $\hat\phi_\eps$ is a normalized Haar functional
(Definition \ref{Haar things}), its centrality follows from
Proposition~\ref{Evaluation}. Then we have for any $x\in\hat M$:
\begin{eqnarray*}
\hat\phi_\eps(xx^*) 
&=& <x\otimes x^*,\Delta(p_{\eps})>
 = \sum_{i\in\calI} \, \frac{1}{d_i}\,
     \sum_{kl} \, e_{kl}^{(i)}(x)\,{S}e_{lk}^{(i)}(x^*) \\
&=& \sum_{i\in\calI} \, \frac{1}{d_i}\,\sum_{kl}\,\left|
    \,e_{kl}^{(i)}(x) \, \right|^2  \geq  0,
\end{eqnarray*}
which equals to $0$ iff $x=0$.

The second statement now is clear from Remark \ref{selfduality of axioms}a) since
$\hat\phi_\eps$ gives a faithful representation of $\hat M$ which means that
$\hat M$ is a $C^*$-algebra.

\begin{corollary}
\label{Faithful trace} Given a weak Kac algebra, there exists a unique
normalized Haar trace $\phi_\eps = p_{\hat\eps}$.
\end{corollary}

{\em Proof.} The existence is clear from the previous proof, the
uniqueness follows from Theorem \ref{Haar Theorem}.

\begin{remark} Let us describe all Haar traces on a weak Kac algebra.
Proposition~\ref{Evaluation} shows that
$\varsigma\Delta(p_i)=\Delta(p_i)\ (\forall i\in \calS)$, therefore,
any linear combination $\sum_{i\in\calS}\lambda_ip_i$ with positive
coefficients $\lambda_i$ defines a trace on $\hat M$.
But the proof of Theorem~\ref{Haar Theorem} shows that any Haar trace
is of this form. Such a trace is faithful iff all $\lambda_i > 0$.
\end{remark}

\end{subsection}

\begin{subsection} {Haar conditional expectations}

\begin{proposition}
Given a weak Kac algebra $(M,\Delta,S,\varepsilon)$ and a normalized
Haar
trace $\phi_\eps$ on it, there exist unique faithful target and source
Haar
conditional expectations $E_t:M\to N_t$ and $E_s:M\to N_s$ such that
$\phi_\eps=\phi_\eps\circ E_t$ and $\phi_\eps=\phi_\eps\circ E_s$.
They could be defined as follows:
$$E_t(x)=(id\otimes\phi_\eps)\Delta(x)=(S\otimes\phi_\eps)((1\otimes
x)e),
E_s(x)=(\phi_\eps\otimes id)\Delta(x)\ (\forall x\in M).$$
We also have $E_t\circ S=S\circ E_s$ and
$$(E_t\otimes E_t)((\Delta(x)(y\otimes z))=\varsigma(E_t\otimes
E_t)((S(y)\otimes x)\Delta(z))\ (\forall x,y,z\in M).$$
\end{proposition}

{\em Proof.}
Clearly both $E_t$ and $E_s$ are linear unital $*$-maps from $M$ to
$N_t$ and
$N_s$ respectively such that $E_t|N_t=id_{N_t}, E_s|N_s= id_{N_s}$, and
$E_t\circ S=S\circ E_s$. We also have $\forall x\in M,n\in N_t$:
\begin{eqnarray*}
E_t(x)\cdot n
&=& (S\otimes\phi_\eps)((1\otimes x)e)\cdot n
=  (S\otimes\phi_\eps)((1\otimes x)(S(n)\otimes 1)e) \\
&=& (S\otimes\phi_\eps)((1\otimes xn)e) = E_t(xn)
\end{eqnarray*}
and
$$
\phi_\eps\circ E_t(x)=(\phi_\eps\circ S\otimes\phi_\eps)((1\otimes
x)e)=\phi_\eps(x).
$$
Now it follows from (\cite{GHJ}, 2.6.2) that $E_t$ is indeed a faithful
conditional expectation. The relations $(id\otimes
E_t)\circ\Delta=\Delta
\circ E_t$ and $(E_s\otimes id)\circ\Delta=\Delta\circ E_s$ are obvious
from the coassociativity of $\Delta$.
In order to prove the last relation, it suffices to show that
$(\alpha\otimes\beta)[(E_t\otimes E_t)(\Delta(x)(y\otimes z))]=
(\alpha\otimes\beta)[\varsigma(E_t\otimes E_t)(S(y)\otimes
x)\Delta(z))]$
for any linear functionals $\alpha,\beta$ on $N_t$.
But these functionals can be represented as
$\alpha(\cdot)=\phi_\eps(n\cdot),\ \beta(\cdot)=\phi_\eps(n'\cdot)$ for
some
$n,n'\in N_t$. Then, using the properties of traces, conditional
expectations,
the invariancy of $\phi_\eps$ and Propositions~\ref{Cartan properties},
\ref{crucial}, we have $\forall x,y,z\in M$:
\begin{eqnarray*}
&& (\alpha\otimes\beta)[(E_t\otimes E_t)(\Delta(x)(y\otimes z))] = \\
&=& (\phi_\eps\otimes\phi_\eps)[(n\otimes n')(E_t\otimes
E_t)(\Delta(x)(y\otimes z))] \\
&=& (\phi_\eps\otimes\phi_\eps)[(E_t\otimes E_t)((n\otimes
n')\Delta(x)(y\otimes z))] \\
&=& (\phi_\eps\otimes\phi_\eps)[(n\otimes n')\Delta(x)(y\otimes z)] \\
&=& (\phi_\eps\otimes\phi_\eps)[(n\otimes 1)e\Delta(x)(y\otimes zn')]
\\
&=& (\phi_\eps\otimes\phi_\eps)[\Delta(nx)(y\otimes zn')] \\
&=& (\phi_\eps\otimes\phi_\eps)[(S(y)\otimes nx)\Delta(zn')] \\
&=& (\phi_\eps\otimes\phi_\eps)[(S(y)\otimes nx)\Delta(z)(n'\otimes 1)]
\\
&=& (\phi_\eps\otimes\phi_\eps)[(n'\otimes n)(S(y)\otimes x)\Delta(z)]
\\
&=& (\phi_\eps\otimes\phi_\eps)[(E_t\otimes E_t)((n\otimes n')
    \varsigma((S(y)\otimes x)\Delta(z))]  \\
&=& (\phi_\eps\otimes\phi_\eps)[(n\otimes n')
    \varsigma((E_t\otimes E_t)((S(y)\otimes x)\Delta(z))] \\
&=& (\alpha\otimes\beta)[\varsigma(E_t\otimes E_t)
    ((S(y)\otimes x)\Delta(z))]
\end{eqnarray*}

\begin{remark}
For a general Haar trace $\phi$ the conditional expectations
$E_t$ and $E_s$ satisfying all the above properties are
the orthogonal projectors onto $N_t$ and $N_s$ respectively in the Hilbert
space given by the GNS-construction for $(M,\phi)$, but they cannot be
given by the above formulae (see \cite{Val3}).
\end{remark}

\begin{proposition}
\label{positivity}
a) In a weak Kac algebra any of the relations $e(1\otimes x)=0,\
e(x\otimes 1)=0$ implies $x=0$.

b) If $y\in M$ and one of the elements $e(1\otimes y),\ e(y\otimes 1),
\ (1\otimes y)e,\ (y\otimes 1)e$ is positive, then $y$ is positive.
\end{proposition}

{\it Proof.}
a) Since $\mu(S\otimes id)(e(1\otimes x))=x\ (\forall x\in M)$, the
first of the statements a) is clear. The second one can be proved
similarly.

b) If $e(y\otimes 1)$ is positive and $\phi_\eps$ is the normalized
Haar trace,
then:
$$
(\phi_\eps\otimes\phi_\eps)((x^*\otimes 1)e(y\otimes 1)(x\otimes 1))=
\phi_\eps(x^*yx)\geq 0\ (\forall x\in M).
$$
This means that $y$ is positive. Now all other statements are clear.

\begin{proposition} (cf. \cite{Val3})
a)  the map $y\to E^o_t(y)=\mu(S\otimes id)((1\otimes y)e)$
is a  faithful conditional expectation from $M$ to $M\cap N'_t$
satisfying
the relations
$$
e(1\otimes y)e=e(1\otimes  E^o_t(y))=(1\otimes  E^o_t(y))e.
$$

b) For any fixed Haar trace $\phi,\ E^o_t(y)$ is the unique faithful
conditional
expectation from $M$ to $M\cap N'_t$ such that $\phi=\phi\circ E^o_t$.
\end{proposition}

{\em Proof.}
a) For any representation $e=\sum_{i=1}^n x_i\otimes y_i$ with minimal
possible $n$ we have, using Proposition~\ref{crucial}:
\begin{eqnarray*}
(1\otimes  E^o_t(y))e
&=& (1\otimes \sum_{i=1}^n S(x_i)y y_i)e
    =\sum_{i=1}^n (1\otimes S (x_i)y) (1\otimes y_i)e \\
&=& \sum_{i=1}^n (S(y_i)\otimes S (x_i)y)e
    =\varsigma (S\otimes S )(e)(1\otimes y)e \\
&=& e(1\otimes y)e\ (\forall y\in M). \\
\end{eqnarray*}
Similarly one can prove the other relation. These relations and
Proposition~\ref{crucial} imply for any $y\in M,\ n\in N_t$:
\begin{eqnarray*}
(1\otimes  E^o_t(y)n)e
&=& (1\otimes  E^o_t(y))(1\otimes n)e
=(S(n)\otimes 1)(1\otimes  E^o_t(y))e \\
&=& (S(n)\otimes 1)e(1\otimes y)e
 = (1\otimes n)e(1\otimes y)e \\
&=& (1\otimes nE^o_t(y))e.
\end{eqnarray*}
Using Proposition~\ref{positivity}~a), we see that
$E^o_t(y)n=nE^o_t(y)$,
so $ E^o_t(y)\in M\cap N'_t$. If $x,x'\in M\cap N'_t, y\in M$, then
$$
E^o_t(xyx')=\sum_{i=1}^n S(x_i)(xyx') y_i
=x(\sum_{i=1}^n S(x_i)y y_i)x'=xE^o_t(y)x'.
$$
If now $y=zz^*\in M$ is positive, then
$$
e(1\otimes  E^o_t(zz^*)) = e(1\otimes zz^*)e =
e(1\otimes z)[e(1\otimes z)]^*\geq 0.
$$
So, according to Proposition~\ref{positivity}~b),
$E^o_t(zz^*))$ is positive, and  if $E^o_t(zz^*))=0$,
then Proposition~\ref{positivity}~a) implies $z=0$.
Since $E^o_t(1)=\varepsilon_t(e)=1$, then $E^o_t$ is indeed the
faithful conditional expectation from $M$ to $M\cap N'_t$.

b) For any $x\in M$ one has:
$\phi\circ E^o_t(x)=\phi(\sum_{i=1}^n y_iS(x_i)x)=\phi(x)$,
so the uniqueness of $E^o_t$ is a consequence of (\cite{GHJ}, 2.6.2.)
\end{subsection}

\begin{subsection}{Generalized Kac algebras}

Let ${\cal K}=(M,\Delta,S,\phi)$ be a gen.\ Kac algebra and
$\hat{\cal K}=(\hat M,\hat\Delta,\hat S,\hat\phi)$ the corresponding
dual gen.\ Kac algebra \cite{Y}. As a linear space, $\hat M$ can be
identified with $M$ because any linear functional $\alpha$ on $M$ can
be represented as $\alpha(x)=\phi(xy)=\phi(yx)$ for some $y\in M.$
Let us consider the following linear functional
$\varepsilon:M\to C:\ \varepsilon(x)=\phi(\hat 1 x)\ (\forall x\in M)$,
where $\hat 1\in M$ corresponds to the unit in $\hat M$. From (\cite{Y},
\S 5) we have $S(\hat 1)=(\hat 1)^*=\hat 1$,
$$
(\phi\otimes\phi)[\Delta(\hat 1)(x\otimes y)]=\phi(xS(y))\
(\forall x,y\in M),
$$
which implies respectively $\varepsilon\circ S(x)=\varepsilon(x),\
\varepsilon(x^*)=\overline{\varepsilon(x)}$ and
$(id\otimes\varepsilon)[\Delta(x)]= x\ (\forall x\in M)$.
The last one gives: $(\varepsilon\otimes id)[\Delta(x)] =
x\ (\forall x\in M)$.

Now let us show that $(M,\Delta,S,\varepsilon)$ is a weak Kac algebra.
According to Proposition~\ref{equivalence of axioms}~, it suffices to
verify axioms A2-A4.

One can calculate, using the invariancy of $\phi$ and the latest
relations:
\begin{eqnarray*}
(\phi\otimes id)[\Delta(\hat 1)(x\otimes y)]
&=& (\phi\otimes id)[\Delta(\hat 1)(x\otimes 1)]y \\
&=& (\phi\otimes S )[(\hat 1\otimes 1)\Delta(x)]y \\
&=& [S\circ(\varepsilon\otimes id)\Delta(x)]y=S(x)y\
(\forall x,y\in M),
\end{eqnarray*}
which gives $(\phi\otimes id)[\Delta(\hat 1)A]=
\mu(S\otimes id)(A)\ (\forall A\in M\otimes M)$.
Similarly we have: $(id\otimes\phi)[A\Delta(\hat 1)] =
\mu(id\otimes S )(A)\ (\forall A\in M\otimes M)$.

Let us first verify A4, using the invariancy of $\phi$ and the latest
relations:
\begin{eqnarray*}
(id\otimes\varepsilon)[(1\otimes x)\Delta(1)]
&=& (id\otimes\phi)[(1\otimes\hat 1x)\Delta(1)]
= (S\otimes\phi)[\Delta(\hat 1x)] \\
&=& \mu(S\otimes id)\Delta(x)\quad (\forall x\in M).
\end{eqnarray*}
Then let us verify A2, using the same reasoning and also A4:
\begin{eqnarray*}
(\varepsilon\otimes id)[(x\otimes 1)\Delta(y)]
&=& (\phi\otimes id)[(\hat 1 x\otimes 1)\Delta(y)] \\
&=& (\phi\otimes S)[\Delta(\hat 1 x)(y\otimes 1)] \\
&=& (\phi\otimes S)[\Delta(\hat 1)\Delta(x)(y\otimes 1)] \\
&=& S[\mu(S\otimes id)(\Delta(x)(y\otimes 1))] \\
&=&  S[\mu(S\otimes id)\Delta(x)]y  \\
&=& (\phi\otimes S)[\Delta(\hat 1 x)]\cdot y \\
&=& (\phi\otimes id)[(\hat 1 x\otimes 1)\Delta(1)]\cdot y \\
&=& (\varepsilon\otimes id)[(x\otimes 1)e(1\otimes y)].
\end{eqnarray*}

Finally, it suffices to verify A3* (see Remark~\ref{selfduality of
axioms}~b)). But this is easy if one uses the following results
obtained in \cite{Val3}: a) $N_t=M\cap\hat M$ is a target Cartan
subalgebra;
b) $e\in N_s\otimes N_t$. Indeed:
\begin{eqnarray*}
(id\otimes\varepsilon\otimes id)[(e\otimes 1)(1\otimes e)]
&=& (id\otimes\varepsilon)(e(1\otimes 1_{(1)}))\otimes 1_{(2)} \\
&=& (id\otimes\varepsilon)(\Delta(1_{(1)}))\otimes 1_{(2)}=e.
\end{eqnarray*}
Thus, we get the following

\begin{theorem}
\label{generalized2weak}
Let $(M,\Delta,S,\phi)$ be a gen.\ Kac algebra,
$\varepsilon(x)=\phi(\hat 1 x)\ (\forall x\in M)$,
where $\hat 1\in M$ is the unit in $\hat M$.
Then $(M,\Delta,S,\varepsilon)$ is a weak Kac algebra.
\end{theorem}

\begin{corollary}
\label{unique normalized}
On any gen.\ Kac algebra there is a  unique Haar normalized
trace.
\end{corollary}

Indeed, to every gen.\ Kac algebra there corresponds a weak Kac algebra
(Theorem~\ref{generalized2weak}) with the unique normalized Haar trace 
(Corollary~\ref{Faithful trace}).

\begin{definition}
A homomorphism $\pi$ of gen.\ Kac algebras with normalized Haar
traces $\psi_1$ and $\psi_2$ is such a homomorphism of their gen.\
coinvolutive $C^*$-bialgebras (Definition~\ref{Def: generalized
Hopf}~c)), that $\psi_2\circ\pi=\psi_1$.
\end{definition}
Clearly gen.\ Kac algebras form a category; Corollary~\ref{unique normalized}
and Theorem~\ref{generalized2weak} describe its equivalence to the category
of weak Kac algebras.

\end{subsection}

\begin{subsection}{Kac bimodules}

The structure of a gen.\ counital $C^*$-bialgebra is not selfdual: 
there is a unit but not necessarily a counit. Let us show that certain extra
condition motivated by Lemma \ref{epsilon is a trace} gives a weak Kac algebra
structure.

Let $\theta_t$ and $\theta_s$ be the traces of left regular representations 
of Cartan subalgebras $N_t$ and $N_s$ respectively. Since $S:N_s\to N_t$ is 
an anti-$*$-isomorphism, we have $\theta_t\circ S=\theta_s$. Also
$\theta_s\circ\eps_s (x)=\theta_s(Sx_{(1)}x_{(2)})=\theta_t\circ S
(Sx_{(1)}x_{(2)})=\theta_t(Sx_{(2)}x_{(1)} )=\theta_t(x_{(1)}Sx_{(2)})
=\theta_t\circ\eps_t(x)\ (\forall x\in M).$

Denoting $\eps:=\theta_s\circ\eps_s=\theta_t\circ\eps_t:M\to \Bbb{C}$, 
we have immediately $\eps\circ S =\eps$ and, 
using Remark~\ref{source counit} d) and self-adjointness of $\theta_t$:
$$ 
\eps(x) = \theta_t\circ(\eps_t(x^*)^*) =
\overline{\theta_t\circ\eps_t(x^*)}=\overline{\eps(x^*)}.
$$

\begin{theorem}
\label{characterization} 
A gen.\ counital $C^*$-bialgebra $(M,\Delta,S)$ is a weak Kac algebra 
if and only if  $(id\otimes\eps)\circ\Delta=id$. 
\end{theorem}

{\it Proof.} Any weak Kac algebra is a gen.\ counital $C^*$-bialgebra 
satisfying the above condition (Corollary~\ref{target and source}, 
Proposition~\ref{epsilon is a trace}). Conversely:

1) Using Remark \ref{source counit} c) and the condition of the
theorem, we have the axiom A3* of a weak Kac algebra:
$$
(id\otimes\eps\otimes id)[(e\otimes 1)(1\otimes
e)]=(id\otimes\eps\otimes id)(\Delta\otimes id)(e)=e.
$$

2) From Remark~\ref{source counit} d) and Proposition \ref{crucial} we have:
$$
(\eps\otimes id)[\Delta(y)(x\otimes 1)]= (\theta_t\circ\eps_t\otimes
id)[\Delta(y)(S\circ\eps_t(x)\otimes 1)]=y\eps_t(x)\ (x,y\in M);
$$
for $y=1$ this gives the axiom A4' of a weak Kac algebra equivalent to
A4.

3) Using the above relation 2) and the axiom A4', we have the axiom A2'
of a weak Kac algebra (which is equivalent to A2):
$$
(\eps\otimes id)[\Delta(y)(x\otimes 1)]=y\eps_t(x)=(\eps\otimes
id)((1\otimes y)e(x\otimes 1)),
$$
The final result follows from Proposition~\ref{equivalence of axioms} 
and Remark~\ref{selfduality of axioms} b).

Let us remark that the condition of Theorem \ref{characterization} makes 
sense also for counital Hopf bimodules (see step a) of the proof of Lemma~\ref{unital equivalence}).

\begin{definition}
\label{Kac bimodules} A counital Hopf bimodule $B=(N,M,t,\Delta,S)$ is said 
to be a {\em Kac bimodule} if it satisfies the condition of 
Theorem~\ref{characterization}.
\end{definition}

Now, taking into account Lemma~\ref{unital equivalence} and Theorem
\ref{characterization}, we get
\begin{theorem} The categories of Kac bimodules and weak Kac algebras
are equivalent.
\end{theorem}
\end{subsection}
\end{section}
\begin{section} {Examples}

\begin{subsection}{Actions of finite groups on weak Kac
algebras}\label{action}
\begin{definition}\label{Def: action}
A right action of a finite group $G$ on a weak Kac algebra 
$WK=(M,\Delta,S,\eps)$ is a right action $M\triangleleft G$ of $G$ on a
$C^*$-algebra $M$ such that, for all $m\in M,\ g\in G$ :
$$
\Delta(m\triangleleft g)=\Delta(m)\triangleleft (g\otimes g),\ S(m\triangleleft g)=S(m)\triangleleft g,\
\eps(m\triangleleft g)=\eps(m).
$$
\end{definition}
This is exactly a right action of a cocommutative Kac algebra
${\cal K}_s(G) =({\rm C}G,\Delta_s,S_s,\varepsilon_s)$ with
$\Delta_s(g)=g\otimes g, S_s(g)=g^{-1},
\varepsilon(g)=1$ \cite{ES} on $M$ such that
$\Delta(m\triangleleft\tilde{g})=\Delta(m)\triangleleft
\Delta_s(\tilde{g}),\ S(m\triangleleft\tilde{g})=S(m)\triangleleft\tilde{g},\
\eps(m\triangleleft\tilde g)=\eps(m)\
(\forall\tilde{g}\in {\rm C}G).$

Given a right action of $G$ on WK, let us consider a collection $WK_G =
(M_G,\Delta_G,S_G,\eps_G)$. Here $M_G=M\#G $ is a {\it crossed product
$C^*$-algebra} \cite{M} generated by elements of the form
$m\otimes g\ (\forall m\in M, g\in G)$ with operations
$$
(m\otimes g)(n\otimes h)=(m\triangleleft h)n\otimes gh,\ 1_{M_G}=1\otimes 1_G,
(m\otimes g)^*=(m\triangleleft g^{-1})^*\otimes g^{-1},
$$
where $m\in M,\ g\in G,\ 1_G$ is a unit of $G$. The {\it coproduct}
$\Delta_G$, {\it antipode} $S_G$, and {\it counit} $\eps_G$ are
defined respectively by: $\Delta_G(m\otimes
g)=(id\otimes\varsigma\otimes id_G) (\Delta(m)\otimes g\otimes g)$,
 $S_G(m\otimes g)=S(m\triangleleft g^{-1})\otimes g^{-1}$, where
$\varsigma:M\otimes
{\rm C}G\to{\rm C}G\otimes M$ is the usual flip, and $\eps_G(m\otimes
g)=\eps(m)$.

\begin{lemma}
\label{Lem: action} $WK_G = (M_G,\Delta_G,S_G,\eps_G)$ is
a weak Kac algebra.
\end{lemma}
{\it Proof.} a) $\Delta_G$ is coassociative as tensor product of two
coassociative coproducts
$\Delta$ and $\Delta_s$ and obviously multiplicative. Let
$m,n\in M, g,h\in G$.

Relations $S_G(1_{M_G})=1_{M_G}$ and $S^2_G=id_{M_G}$ are obvious; we
also have:
$S_G[(m\otimes g)(n\otimes h)]=S_G((m\triangleleft h)n\otimes gh)=S(((m\triangleleft
h)n)\triangleleft(gh)^{-1}\otimes(gh)^{-1}=
S((m\triangleleft g^{-1})(n\triangleleft(gh)^{-1}))\otimes(gh)^{-1}$.

On the other hand :
$S_G(n\otimes h)S_G(m\otimes g)=(S(n\triangleleft h^{-1})\otimes h^{-1})(S(m\triangleleft
g^{-1})\otimes g^{-1})=
([(S(n\triangleleft h^{-1})]\triangleleft g^{-1})(S(m\triangleleft g^{-1}) \otimes h^{-1}g^{-
1})=S((m\triangleleft g^{-1})(n\triangleleft (gh)^{-1}))\otimes (gh)^{-1}$. Thus, $S_G$
is antimultiplicative.

$(\Delta_G\circ S_G)(m\otimes g)=\Delta_G(S(m\triangleleft g^{-1})\otimes g^{-
1})=
(id\otimes\varsigma\otimes id_G)(\varsigma_M (S\otimes S)\Delta(m\triangleleft
g^{-1})\otimes g^{-1}\otimes g^{-1})= \varsigma_{M_G}(S_G\otimes
S_G)\Delta_G(m\otimes g)$,
where the usual flips act: $\varsigma_{M}$ in $M\otimes M$ and
$\varsigma_{M_G}$ in $M_G\otimes M_G$.

Thus, $(M_G,\Delta_G,S_G)$ is a gen.\ coinvolutive $C^*$-bialgebra.

b) Clearly, $\eps_G$ is a counit with respect to $\Delta_G$. Let us verify 
Axiom 2):
\begin{eqnarray*}
\eps_G((m\otimes g)(n\otimes h))
&=& \eps((m\triangleleft h)n)=
(\eps\otimes\eps)[((m\triangleleft h)\otimes 1)e(1\otimes n)] \\
&=&(\eps\otimes\eps)([((m\otimes 1)e)\triangleleft(h\otimes h)](1\otimes n))\\
&=& (\eps_G\otimes\eps_G)[(m\otimes g\otimes
1_{M_G})\Delta_G(1_{M_G})(1_{M_G}\otimes n\otimes h)].
\end{eqnarray*}
c) Finally, let us verify Axiom 3) ($\mu_G$ is the multiplication
in $M_G$):
\begin{eqnarray*}
& & (\mu_G\otimes id_{M_G}))(S_G\otimes id_{M_G}\otimes id_{M_G})
(\Delta_G\otimes id_{M_G})\Delta_G(m\otimes g)  = \\
&=& (\mu_G\otimes id_{M_G})(S_G\otimes id_{M_G}\otimes
id_{M_G})(m_{(1)}\otimes g\otimes m_{(2)}\otimes g
\otimes m_{(3)}\otimes g) \\
&=& (id\otimes\varsigma\otimes id_G)[(\mu\otimes id)(S\otimes id\otimes
id)(\Delta\otimes id)\Delta(m)\otimes 1_G\otimes g] \\
&=& (id\otimes\varsigma\otimes id_G)[(1\otimes m)e\otimes 1_G\otimes g]
= (1\otimes 1_G\otimes m\otimes g)\Delta_G(1\otimes 1_G).
\end{eqnarray*}

\begin{remark} One can show that the Haar trace on $WK_G$ is just the
dual trace for the Haar trace on $WK$ with respect to the action of $G$:
$$\phi_G(m\otimes g)=\delta_{1_G,g}\phi(m)\ (\forall m\in M,g\in G).$$
\end{remark}
Now we describe a series of concrete examples of nontrivial (non-commu-
tative and non-symmetric) weak Kac algebras of dimensions $n^3\ (n\geq 3)$.

Let us start with a weak Kac algebra $M=\Bbb{C}(K_n)$ of functions on the
transitive principal groupoid $K_n$ on a set of $n$ elements (for various
examples of groupoids see \cite{R}). 
Then $M=\mbox{span}\{e_{ij}\}^n_{i,j=1}$ is
a commutative $C^*$-algebra of functions on the set $\{e_{ij}\}$ equipped
with the coproduct
$\Delta:e_{ij}\to\sum^n_{k=1}e_{ik}\otimes e_{kj}$, the antipode
$S:e_{ij}\to e_{ji}$ and the
counit $\eps:e_{ij}\to\delta_{ij}$ ($i,j,k\in\{1,...,n\},\ \delta_{ij}$
is the Kronecker symbol).

Let us consider the action of the cyclic group $G={\rm Z}/n{\rm Z}$ on this
weak Kac algebra
given by $e_{ij}\triangleleft\alpha:=e_{i+1,j+1}\ (\alpha$ is a generator of
$G,\ i,j\in\{1,...,n\}$ and the summation is modulo $n$). Applying Lemma
\ref{Lem: action} (all its conditions are satisfied), we get a description
of the resulting weak Kac algebra in terms of matrix units
$f^k_{ij}:=e_{j,j+k}\otimes\alpha^{j-i}$ of $C^*$-algebra $M_G=M\#
G\cong\oplus^n_{k=1}\,M_n(\Bbb{C})$:
$$\Delta_G(f^k_{ij})=\sum^n_{r=1}f^r_{ij}\otimes f^{k-r}_{i+r,j+r},\
S_G(f^k_{ij})=f^{n-k}_{i+k,j+k},\
\eps_G(f^k_{ij})=\delta_{n,k},$$
($i,j,k\in\{1,...,n\}$ and the summation is modulo $n$. Then one can show
that:

a) the above weak Kac algebra is non-trivial;

b) its Haar trace is the canonical trace on $M_G$;

c) $\eps_t(f^k_{ij})=\delta_{n,k}\sum^n_{r=1}f^{n-r}_{j+r,j+r},\
\eps_s(f^k_{ij})=\delta_{n,k}\sum^n_{r=1}f^{r}_{i,i}$;

d) $N_t = \mbox{span}\{\sum^n_{r=1}f^{n-r}_{j+r,j+r}\},\
N_s = \mbox{span}\{\sum^n_{r=1}f^{r}_{i,i}\}$.

For $n=2$ this non-trivial weak Kac algebra has dimension $8$. For a
non-trivial Kac algebra of dimension $8$ see \cite{KP}. Let us show that this
dimension is minimal possible.

\begin{proposition}
All weak Kac algebras of dimension $<8$ are either commutative or
cocommutative (therefore, come from groupoids or groups).
\end{proposition}

{\em Proof.} This statement is well-known for usual Kac algebras, so we
may assume that dimension of Cartan subalgebras $\geq 2$. Let $M$ be a
non-commutative weak Kac algebra of dimension $<8$, then $M\cong\Bbb{C}^k
\oplus M_2(\Bbb{C}),\, 0\leq k\leq 3$ as a $C^*$-algebra. Let us show
that $M$ is cocommutative.

For $k=0$ this follows from the classification of elementary weak Kac
algebras in section 3.2. For $k>0$ it suffices to prove that $M$ is a
direct sum of two weak Kac algebras, since any weak Kac algebra of
dimension $\leq 3$ is cocommutative. Clearly, the last property is
equivalent to existence of a non-trivial projection
$q$ in $N_s\cap N_t \cap \mbox{Center}(M)$, 
the hyper-center of $M$ \cite{NSzW}.

Due to Proposition~\ref{Evaluation} and the multiplicativity of
$\Delta$, for the counital quotient $M_\eps \subset M$
(Proposition~\ref{Canonical quotient}) we have
$M_\eps \neq M_2(\Bbb{C})$, and $M_\eps \neq \Bbb{C}$ since
it must contain a copy of a Cartan subalgebra by Remark 2.1.12. Using
Proposition~\ref{Evaluation} again, we see that all
other possibilities lead to the existence of a proper minimal
subprojection $p_i$ (Remark~\ref{support}) of the Haar projection
$p_\eps\in M$, such that $\varsigma\Delta(p_i) = \Delta(p_i)$ (i.e.
$p_i$ is cocommutative). Since $S(p_i) = p_i$, Proposition~\ref{mult}
implies  $(S\otimes\id)\Delta(p_i) = (\id\otimes S)\Delta(p_i)$ and
$q =\eps_s(p_i) = \eps_t(p_i)$ is a non-trivial projection in the
hyper-center of $M$.
\medskip
\end{subsection}

\begin{subsection}{Elementary weak Kac algebras}

In Subsection 2.2 we have seen that every counital quotient of a
weak Kac algebra is a direct sum of weak Kac algebras which are
simple as algebras.

\begin{define}
\label{elementary}
A weak Kac algebra is said to be {\em elementary} if it is
simple as an algebra, i.e., is isomorphic to the full
matrix algebra $M_n(\Bbb{C})$.
\end{define}
\medskip

The next theorem classifies all elementary weak Kac algebras and shows
that any f-d $C^*$-algebra can appear as a Cartan
subalgebra of an elementary weak Kac algebra in a unique way.

\medskip

%
\begin{theorem}
\label{Classification of elementary weak Kac algebras}
Given a f-d $C^*$-algebra $A$, there exists a unique elementary weak Kac
algebra $M(A) = (M_n(\Bbb{C}),\Delta,S,\eps)\ (n=\dim A)$ with Cartan
subalgebras isomorphic to $A$.
\end{theorem}

{\em Proof.} 1) For $A=\oplus_\alpha\, M_{n_\alpha},\
n=\sum_\alpha\,n_\alpha^2$
let us write down $M(A)$ explicitly, if matrix units
$\{E_{ij\alpha}^{kl\beta} \}_{i,j=1\dots n_\alpha}$,
$$
E_{ij\alpha}^{kl\beta}\,E_{\tilde{k}\tilde{l}\tilde{\beta}}^{\tilde{i}
\tilde{j}\tilde{\alpha}}
=
\delta_{\beta\tilde{\beta}}\,\delta_{k\tilde{k}}\,\delta_{l\tilde{l}}\,
E_{ij\alpha}^{\tilde{i}\tilde{j}\tilde{\alpha}},
\qquad
(E_{ij\alpha}^{kl\beta})^* =E_{kl\beta}^{ij\alpha}.
$$
Equip $M_n(\Bbb{C})$ with the following structure of a weak
Kac algebra :
\begin{eqnarray*}
\Delta(E_{ij\alpha}^{kl\beta})
&=& \frac{1}{\sqrt{n_\alpha n_\beta} }
\,\sum_{uv} E_{iu\alpha}^{kv\beta} \otimes E_{uj\alpha}^{vl\beta},\\
 S(E_{ij\alpha}^{kl\beta})
&=& E_{lk\beta}^{ji\alpha}, \\
\varepsilon(E_{ij\alpha}^{kl\beta})
&=&  \sqrt{n_\alpha n_\beta}\,\delta_{ij}\,\delta_{kl}.
\end{eqnarray*}
The verification of axioms is a straightforward computation.
We have
$$
e = \Delta(1) =
\sum_{\alpha}\, \frac{1}{n_\alpha}\, \sum_{pq}\,
(\sum_i E_{ip\alpha}^{iq\alpha}) \otimes (\sum_j
E_{pj\alpha}^{qj\alpha}),
$$
therefore, $F_{pq \alpha} = \sum_i\, E_{ip\alpha}^{iq\alpha}$ and
$G_{pq\alpha} = \sum_j\, E_{qj\alpha}^{pj\alpha}$ are matrix units
for the source and target Cartan subalgebras respectively, isomorphic to
$A$.

2) The uniqueness. If $\tilde{M}(A) = (M_n(\Bbb{C}),\,\tilde\Delta,\,
\tilde S,\,\tilde\eps)$
is another weak Kac algebra with $\tilde{N_s}\cong\tilde{N_t}
\cong \oplus_\alpha\, M_{n_\alpha}(\Bbb{C})$, then, by Remark~\ref{support},
$\sum\nolimits_{\alpha=1}^{K}\, n^2_\alpha = \dim A =\deg\pi_\eps = n$.

The inclusion $N_s\subset M$ is determined up to conjugation
by the vector $(d_1,\dots d_K)$ where $d_\alpha$ is the multiplicity of
$M_{n_\alpha}(\Bbb{C})$ in $M$.  By Proposition \ref{Cartan
subalgebras},
$N_t \subset N_s'\cap M = \oplus_\alpha\, M_{d_\alpha}(\Bbb{C})$,
therefore we may assume that $n_\alpha \leq d_\alpha$. Since the
inclusion
$N_s\subset M$ is unital,  $\sum d_\alpha n_\alpha =n$, so
$0 = \sum d_\alpha n_\alpha - \sum n_\alpha^2 =
\sum (d_\alpha - n_\alpha)n_\alpha$,
and $d_\alpha = n_\alpha\ (\forall\alpha)$.
Thus,  $N_t = N_s'\cap M$, and $N_s = N_t'\cap M$.

Hence, we can write $N_s=\oplus_\alpha N_{s\alpha}$ and
$N_t=\oplus_\alpha N_{t\alpha}$, where $N_{s\alpha} \cong
N_{t\alpha} \cong M_{n_\alpha}(\Bbb{C})$, $N_{s\alpha}N_{t\beta} =0$
if $\alpha\neq\beta$, and  $N_{s\alpha}N_{t\alpha} \cong
M_{n^2_\alpha}(\Bbb{C})$. Let $\{F_{pq\alpha}\}_{p,q=1\dots n_\alpha}$
be a system of matrix units in $N_{s\alpha}$, $\alpha =1\dots K$.
Proposition~\ref{formula for e} gives
$$
\tilde{e} =
\sum\nolimits_\alpha\,\frac{1}{n_\alpha}\sum\nolimits_{pq}\,
F_{pq\alpha} \otimes S(F_{qp\alpha}).
$$
Note that 
$ E_{ij\alpha}^{kl\alpha} = \tilde{S}(F_{ki\alpha}) F_{jl\alpha}$ 
are matrix units in
$N_{s\alpha}N_{t\alpha}$ such that
$$
F_{pq\alpha} = \sum\nolimits_i\, E_{ip\alpha}^{iq\alpha},  \qquad
\tilde{S}(F_{qp\alpha}) =  \sum\nolimits_j\, E_{pj\alpha}^{qj\alpha}
$$
Let $\{ E_{ij\alpha}^{kl\beta} \}$ be the system of matrix units
in $M_n(\Bbb{C})$ defined as in 1) which extends the union of the above
systems $\{ E_{ij\alpha}^{kl\alpha} \}$,  $\alpha =1,\dots, K$.
Using the definition of a Cartan subalgebra we compute
\begin{eqnarray*}
\tilde\Delta(F_{pq\alpha})
&=& \frac{1}{n_\alpha}\,
    \sum_{\tilde{p}\tilde{q}}\,
    F_{\tilde{p}\tilde{q}\alpha}\otimes
    E_{\tilde{p}p\alpha}^{\tilde{q}q\alpha}, \\
\tilde\Delta(E_{lp\alpha}^{mq\alpha})
&=& \tilde\Delta(\tilde{S}(F_{ml\alpha})) \tilde\Delta(F_{pq\alpha})
    =  \frac{1}{n_\alpha}\, \sum_{uv}\,
    E_{lu\alpha}^{mv\alpha} \otimes E_{up\alpha}^{vq\alpha}
    = \Delta(E_{lp\alpha}^{mq\alpha}).
\end{eqnarray*}
Since $\tilde\Delta$ is a coassociative $*$-homomorphism, we get
$$
\tilde\Delta(E_{ij\alpha}^{kl\beta}) =
\frac{\lambda_\alpha^\beta}{\sqrt{n_\alpha n_\beta} }
\,\sum_{uv} E_{iu\alpha}^{kv\beta} \otimes E_{uj\alpha}^{vl\beta},
= \lambda_\alpha^\beta \Delta(E_{ij\alpha}^{kl\beta})
$$
for some $\{\lambda_\alpha^\beta \}_{\alpha,\beta =1\dots K}$
with $\vert \lambda_\alpha^\beta \vert =1$,
$ \lambda_\beta^\alpha = \overline{\lambda_\alpha^\beta}$,
and $ \lambda_\alpha^\beta\,\lambda_\beta^\gamma =
\lambda_\alpha^\gamma
\ (\forall \alpha,\,\beta,\,\gamma)$.
The properties of antipode and counit imply:
$\tilde{S}(E_{ij\alpha}^{kl\beta}) =
{(\lambda_\alpha^\beta)}^{-2}\,S(E_{ij\alpha}^{kl\beta})$
and $\tilde\eps(E_{ij\alpha}^{kl\beta}) = {(\lambda_\alpha^\beta)}^{-1}
\eps(E_{ij\alpha}^{kl\beta})$, so the map
$ E_{ij\alpha}^{kl\beta} \mapsto \lambda_\alpha^\beta
E_{ij\alpha}^{kl\beta}$
defines an isomorphism of weak Kac algebras
$M(A) = (M_n(\Bbb{C}),\,\Delta,\,S,\,\eps)$ and
$\tilde{M}(A) = (M_n(\Bbb{C}),\,\tilde\Delta,\,\tilde S,\,\tilde\eps)$.
\medskip

Let us mention that $M(A)$ can be viewed as a quantum transitive principal
groupoid on $A$, since when $A\cong \Bbb{C}^n$ (i.e., when $A$ is abelian)  
$M(A)$ is precisely the groupoid weak Kac algebra of the transitive             
principal groupoid on a set of $n$ elements.    

\begin{corollary}
The dual weak Kac algebra $\hat{M}(A)$ has an algebra structure
$\oplus_{\alpha\beta}\, M_{n_\alpha n_\beta}(\Bbb{C})$; there
are comatrix units
$\{ C_{ij\alpha}^{kl\beta} \}_{i,j=1\dots n_\alpha}$,
$$
\Delta(C_{ij\alpha}^{kl\beta}) = \sum_{pq\gamma}\, 
C_{ij\alpha}^{pq\gamma}\otimes C_{pq\gamma}^{kl\beta},
\qquad
\eps(C_{ij\alpha}^{kl\beta}) =
\delta_{\alpha\beta}\delta_{ik}\delta_{jl},
$$
such that
\begin{eqnarray*}
C_{iu\alpha}^{kv\beta} \cdot
C_{ \tilde{u} j \tilde{\alpha} }^{ \tilde{v} l \tilde{\beta} }
&=& \frac{1}{ \sqrt{n_\alpha n_\beta} }\,\delta_{ \alpha \tilde{\alpha} }
   \delta_{ \beta \tilde{\beta} } \delta_{u\tilde{u}} \delta_{v\tilde{v} }\,
C_{ij\alpha}^{kl\beta},\\
1 &=&\sum_{\alpha\beta\, ij}\,\sqrt{n_\alpha n_\beta}\,
C_{ii\alpha}^{kk\beta},\\
\end{eqnarray*}
$$
{(C_{ij\alpha}^{kl\beta}) }^* = C_{ji\alpha}^{lk\beta},\
S(C_{ij\alpha}^{kl\beta}) = C_{lk\beta}^{ji\alpha}.
$$

\end{corollary}

{\em Proof}. A straightforward computation using the definition
of dual weak Kac algebra.

\end{subsection}
\end{section}
 


\begin{thebibliography}{10}
\bibitem{BS}
S.~Baaj et G.~Skandalis.
\newblock  Unitaires multiplicatifs et dualit\'e pour les produits
crois\'es de $C^*$-alg\`ebres,
\newblock {\em Ann. Sci. Ecole Norm. Sup.}, {\bf 26} (1993), 425--488.

\bibitem{BNSz}
G.~B\"ohm, F.~Nill and K.~Szlach\'anyi.
\newblock  Weak Hopf algebras I: Integral theory and $C^*$-structure,
\newblock {\em Preprint, math.QA}/{\bf 9805116} (1998).

\bibitem{BSz}
G.~B\"ohm, K.~Szlach\'anyi.
\newblock  Weak $C^*$-Hoph algebras: The coassociative symmetry with non-
integral dimensions,
\newblock {\em Banach Center Publications. Warszawa}, {\bf 40} (1997),
9--19;
\newblock A coassociative $C^*$-quantum group with non-integral
dimensions,
\newblock {\em Lett. Math. Phys.}, {\bf 35} (1996), 437--448.

\bibitem{D}
M.-C.~David.
\newblock Paragroupe d'Adrian Ocneanu et algebre de Kac,
\newblock {\em Pacif. J. of Math.}, {\bf 172}, no. 2 (1996), 331--363.

\bibitem{E}
M.~Enock.
\newblock Inclusions irr\'eductibles de facteurs et unitaires
multiplicatifs II,
\newblock {\em J. Funct. Analysis}, {\bf 154}, no. 1 (1998),
67--109.

\bibitem{EN}
M.~Enock and R.~Nest.
\newblock Inclusions of factors, multiplicative unitaries
and Kac algebras,
\newblock {\em J. Funct. Analysis}, {\bf 137} (1996), 466--543.

\bibitem{ES}
M.~Enock and J.-M.~Schwartz.
\newblock {\em Kac algebras and duality for locally compact groups},
\newblock Springer-Verlag, 1989.

\bibitem{EVal}
M.~Enock, J.-M.~Vallin.
\newblock Inclusions of von Neumann algebras and quantum groupoids,
\newblock {\em Inst. Math. de Jussieu}, {Pr\'epublication 156} (1998).

\bibitem{GHJ}
F.M.~Goodman, P. de la Harpe and V.F.R.~Jones.
\newblock Coxeter graphs and towers of algebras, MSRI Publ. 14,
\newblock Springer-Verlag, 1989.

\bibitem{HO}
R.H.~Herman and A.~Ocneanu.
\newblock Index theory and Galois theory for infinite index inclusions
of factors,
\newblock {\em C. R. Acad. Sci. Paris Ser. I Math.}, {\bf 309} (1989),
923--927.

\bibitem{KP}
G.I.~Kac and V.G.~Paljutkin.
\newblock Finite ring groups,
\newblock {\em Trans. Moscow. Math. Society}, {\bf 15} (1966),
251--294.

\bibitem{L}
R.~Longo.
\newblock A duality for Hopf algebras and subfactors I,
\newblock {\em Comm. Math. Phys.}, {\bf 159} (1994), 133--150.

\bibitem{Lu}
J-H.~Lu.
\newblock Hopf algebroids and quantum groupoids,
\newblock {\em Int. J. Math.}, {\bf 7} (1996), no. 1, 47--70.

\bibitem{M}
S.~Majid.
\newblock Foundations of quantum group theory,
\newblock Cambridge University Press, 1995.

\bibitem{Mal}
G.~Maltsiniotis.
\newblock Groupo\"{\i}des quantiques,
\newblock {\em C.R. Acad. Sci. Paris S\'er. I}, {\bf 314} (1992),
249--252.

\bibitem{N}
D.~Nikshych.
\newblock $K_0$-rings and twisting of finite dimensional semisimple
Hopf algebras,
\newblock {\em Communications in Algebra} {\bf 26}, no.1 (1998), 321--
342.

\bibitem{N2}
D.~Nikshych.
\newblock Duality for actions of weak Kac algebras and crossed product
inclusions of II${}_1$ factors,
\newblock {\em Preprint, math.QA}/{\bf 9810049} (1998).

\bibitem{NV}
D.~Nikshych and L.~Vainerman.
\newblock A characterization of depth 2 subfactors of II${}_1$ factors,
\newblock {\em Preprint, math.QA}/{\bf 9810028} (1998).

\bibitem{Nill}
F.~Nill.
\newblock Axioms for weak bilgebras,
\newblock {\em Preprint, math.QA}/{\bf 9805104} (1998).

\bibitem{NSzW}
F.~Nill, K. Szlach\'anyi, H.-W.~Wiesbrock.
\newblock Weak Hopf algebras and reducible Jones inclusions of depth 2.
I: From crossed products to Jones towers,
\newblock {\em Preprint, math.QA}/{\bf 9806130} (1998).
 
\bibitem{O}
A.~Ocneanu.
\newblock A Galois theory for operator algebras (1986).
\newblock {\em Notes} (1986).

\bibitem{R}
J.~Renault.
\newblock A groupoid approach to $C^*$-algebras,
\newblock{\em Lecture Notes in Math.} {\bf 793},
\newblock Springer-Verlag, 1980.

\bibitem{S}
J.-L.~Sauvageot.
\newblock Sur le produit tensoriel relatif d'espaces de Hilbert,
\newblock {\em J. Operator Theory}, {\bf 9} (1983), 237--258.

\bibitem{Szym}
W.~Szymanski.
\newblock Finite index subfactors and Hopf algebra crossed products,
\newblock {\em Proceedings of the AMS}, {\bf 130}, no. 2 (1994),
519--528.

\bibitem{V2}
L.~Vainerman.
\newblock A note on quantum groupoids,
\newblock {\em C.R. Acad. Sci. Paris S\'er. I}, {\bf 315} (1992),
1125--1130.

\bibitem{Val1}
J.-M.~Vallin.
\newblock Bimodules de Hopf et poids op\'eratoriels de Haar,
\newblock {\em J. Operator Theory}, {\bf 35} (1996), 39--65.

\bibitem{Val2}
J.-M.~Vallin.
\newblock Unitaire pseudo-multiplicatif associ\'e \`a un groupo\"{\i}de;
applications \`a la moyennabilit\'e, to be published in
\newblock {\em J. Funct. Analysis}.

\bibitem{Val3}
J.-M.~Vallin.
\newblock Groupo\"{\i}des quantiques finis, in preparation.

\bibitem{Y}
T.~Yamanouchi.
\newblock Duality for generalized Kac algebras and a characterization
of finite groupoid algebras,
\newblock {\em J. Algebra}, {\bf 163} (1994), 9--50.
\end{thebibliography}
\end{document}